\documentclass[11pt,a4paper]{article}

\usepackage[english]{babel}
\usepackage{amsmath}
\usepackage{amssymb}
\usepackage{graphicx}  

\usepackage{epstopdf}       

\usepackage{amsfonts}
\usepackage{color}
\usepackage{verbatim}
\usepackage{caption}
\usepackage{subcaption}
\usepackage{fullpage}

\usepackage{array}

\usepackage{epsfig}
\usepackage{float}

\usepackage{algpseudocode}
\usepackage{url}

\usepackage{theorem}

\usepackage{fullpage}

\usepackage{hyperref}

\usepackage{tikz}

\newtheorem{rem}{Remark}[section]

\newcommand{\EndProofMarker}{$\Box$}

\def\XXint#1#2#3{{\setbox0=\hbox{$#1{#2#3}{\int}$ }
\vcenter{\hbox{$#2#3$ }}\kern-.6\wd0}}

\graphicspath{{./Imgs/}}

\begin{document}
  
\title{An image--informed Cahn--Hilliard Keller--Segel multiphase field model for tumor growth with angiogenesis}
\author{A. Agosti$^{\sharp}$\footnote{Corresponding author. E-mail address: {\tt abramo.agosti@unipv.it} \newline \textit{Email addresses:} {\tt abramo.agosti@unipv.it} (A. Agosti), {\tt alicelucifero@gmail.com} (A. Giotta Lucifero), {\tt sabino.luzzi@unipv.it} (S. Luzzi)}, A. Giotta Lucifero$^{\S}$, S. Luzzi$^{\S,\ddag}$}

\maketitle 

\begin{center}
{\small $^\sharp$  Department of Mathematics, University of Pavia, via Ferrata, 5 - 27100 Pavia, Italy.}
\end{center}

\begin{center}
{\small $^\S$  Neurosurgery Unit, Department of Clinical-Surgical, Diagnostic and Pediatric Sciences, University of Pavia, Viale Brambilla, 74 - 27100 Pavia, Italy.
      }
\end{center}

\begin{center}
{\small $^\ddag$  Neurosurgery Unit, Department of Surgical Sciences, Fondazione IRCCS Policlinico San Matteo, Viale Camillo Golgi, 19 - 27100 Pavia, Italy}
\end{center}

\date{}

\begin{abstract}
In this paper we develop a new four--phase tumor growth model with angiogenesis, derived from a diffuse--interface mixture model composed by a viable tumor component, a necrotic component, a liquid component and an angiogenetic component, coupled with two massless chemicals representing a perfectly diluted nutrient and an angiogenetic factor. This model is derived from variational principles complying with the second law of thermodynamics in isothermal situations, starting from biological constitutive assumptions on the tumor cells adhesion properties and on the infiltrative mechanics of tumor--induced vasculature in the tumor tissues, and takes the form of a coupled degenerate Cahn--Hilliard Keller--Segel system for the mixture components with reaction diffusion equations for the chemicals. The model is informed by neuroimaging data, which give informations about the patient--specific brain geometry and tissues microstructure, the distribution of the different tumor components, the white matter fiber orientations and  the vasculature density. We describe specific and robust preprocessing steps to extract quantitative informations from the neuroimaging data and to construct a computational platform to solve the model on a patient--specific basis. We further introduce a finite element approximation of the model equations which preserve the qualitative properties of the continuous solutions. Finally, we show simulation results for the patient--specific tumor evolution of a patient affected by GlioBlastoma Multiforme, considering two different test cases before surgery, corresponding to situations with high or low nutrient supply inside the tumor, and a test case after surgery. We show that our model correctly predicts the overall extension of the tumor distribution and the intensity of the angiogenetic process, paving the way for assisting the clinicians in properly assessing the therapy outcomes and in designing optimal patient--specific therapeutic schedules. 

\end{abstract}

\vspace*{0.5cm}

\section{Introduction}
The evolution of solid tumors can be separated into an early avascular phase and a vascular phase \cite{intro1}. In the avascular phase, tumor growth is supported by the healthy tissue vasculature and by the diffusion of nutrients through the extravascular space in the peritumoral area. While the size of the tumor increases, the distance between the inner tumor cells and the neighboring blood vessels may exceed the oxygen diffusion distance, causing the cells to become hypoxic. As a consequence of hypoxia conditions, the death rate of tumor cells by apoptosis and necrosis prevail over the mitotic rate, and, in order to provide the necessary infrastructure to further sustain cells proliferation and expansion, the hypoxic tumor cells may gain the ability to release diffusible chemical compounds, known as tumour angiogenesis factors (TAF), to induce the formation of new blood vessels from the surrounding tissues \cite{intro2,intro3}. In response to the angiogenetic stimulus, endothelial cells in the neighbouring normal capillaries are activated to secrete enzymes which degrade their basement membrane, allowing their migration in the extracellular matrix \cite{intro3}. The accumulation of endothelial cells from the parent vessels generates capillary sprouts which grow in length, with the endothelial cells redistributing among sprouts and proliferating at some distance from the sprouts tips. After an initial locally parallel growth of the sprouts, independently from one another, towards the direction of increasing angiogenic signal (via chemotaxis) and along the preferential directions of the basement fibers in the extracellular matrix (via haptotaxis), neighboring sprouts tend to intersecate forming tip--tip and tip--branch anastomoses, giving rise to a new vasculature network with blood circulation which collectively migrates towards the tumor, eventually penetrating into its mass and providing it with an adequate blood supply and microcirculation \cite{intro4,intro5}. 

In the vascular phase of tumor progression, the tumor enhances its ability to invade the surrounding tissues and to metastasise to distant parts of the body.
Tumor--associated microvasculature, being responsible for the delivery of cell substrates (e.g. oxygen and nutrients), metabolite removal and drug delivery, influences in a complex way the response of solid tumors to both systemic and radiation therapy \cite{intro6,intro7,intro9}. Indeed, the highly disorganized structure of tumor--induced vasculature yields an heterogeneous distribution of therapeutics throughout the tumor tissues. Moreover, the efficacy of radiotherapy is highly dependent on tissue oxygenation, which influences the sensitivity of tumor cells to radiotherapy \cite{intro8}, hence it is dependent on the structure and functioning of the tumor vasculature. Radiotherapy itself promotes the angiogenic activity, thus improving in a feedback mechanism the sensitivity of hypoxic cells to future radiotherapy doses \cite{intro10}. Finally, the application of antiangiogenic therapies may have both positive effects, reducing the amount of nutrients reaching the tumor, and negative effects, hindering the delivery of the therapeutics  to the tumor and also increasing the tumor cells motility exacerbating hypoxic conditions. Therefore, understanding the dynamics of the tumor--induced vasculature plays a crucial role in evaluating and possibly designing optimal treatments. In particular, to improve patient outcomes, the therapeutic schedules should be optimized by considering the patient--specific structural and functional characteristics of both the tumor microenvironment and its associated vasculature.

To this end, it seems of crucial importance to study mathematical and computational models, designed on a patient--specific basis, for tumor progression and tumor--induced angiogenesis \cite{intro9}, especially for what concerns those cases that are particularly difficult to be treated with standard therapeutic protocols, such as Glioblastoma Multiforme (GBM). GBM is the most common malignant primary brain tumor in adults. The mean progression-free survival is just over $6$ months. The standard therapy schedule for GBM consists in surgical resection, followed by periodic cycles of radiotherapy and adjuvant chemotherapy \cite{intro11}. Due to the fact that GBM cells infiltrate in the surrounding tissues, following the white matter fiber tracts and the physical structures in the extracellular environment \cite{intro11}, a complete removal of the tumor mass is generally impossible. The treatment with the standard therapy is invariably followed by tumor recurrence, with a mean survival of $14.6$ months. The negative outcome of GBM treatment is also due to its highly vascular nature, determined by endothelial proliferation and the release of angiogenetic factors from the tumor cells. Also, blood-vessel formation in GBM is enhanced through a process called vascular mimicry, where a subpopulation of glioma stem cells differentiates into endothelial cells and participates to the formation of new vessels \cite{intro11}. Indeed, the administration of angiogenesis inhibitors during the therapy schedule has shown to lead to an increase in the progression-free survival.

Treatment evaluation of patients with GBM is important to optimize clinical decisions, and so the treatment--related effects on tumor progression must be addressed and correctly identified. Indeed, clinicians often need to distinguish between progressive disease and pseudoprogression, resulting from neovascularization and from disruption of the blood--brain barrier after chemotherapy. Conventional Magnetic Resonance Imaging (MRI) is unable to distinguish between these two processes, which appear as new enhancement of signal intensity \cite{intro12} in post--contrast T1--weighted images. Moreover, a decrease in enhancement during antiangiogenic treatment may be due to a pseudoresponse of the tumor--related vasculature, and not to a tumor regression induced by cell death \cite{intro12}. For these reasons, a computational model which predict the tumor progression and angiogenetic process under therapy may help in assessing correctly the therapy outcomes and in optimizing the clinical decisions. We note also that recent studies suggest that properly designed radiomics analysis based on Machine Learning tools applied to conventional MRI may help in distinguishing between treatment--related responses and tumor progression \cite{intro13}.     

In light of these observations, biomathematical theories that accurately describe the multiscale dynamic relationships between tumor growth, vasculature and therapy are needed to design optimized therapies and improve the patients outcome. Despite the numerous theoretical and numerical studies in literature which characterize tumor growth and tumor angiogenesis from the cell to the tissue scale, only a small amount of them relate the models with in--vivo imaging data, especially for what concerns the coupling between tumor growth with angiogenesis \cite{intro9}. 

Mathematical models of tumor angiogenesis can be divided into discrete (with endothelial cells represented as individual objects), continuous (with endothelial cells represented as concentrations), or hybrid (combining both the discrete and the continuum pictures) models \cite{intro9}. In discrete models, the evolution in time of the vasculature is dictated by a set of rules governing the individual cells behavior. These models are capable of tracking individual cells and vessel segments. For what concerns continuous models, they describe the evolution of the concentration of endothelial cells over space and time by phenomenological or continuum mechanics laws, describing the overall vessel morphology at the tissue level. These models must be coupled with models for tumor growth and generation of angiogenic stimuli. 

A first attempt to formulate a continuum model of endothelial cells migration following the gradient of TAF (chemotaxis) was introduced in \cite{intro14}, followed by different continuum models generalizations which consider the effects of fibronectin degradation (haptotaxis) and enzymes production \cite{intro15,intro16,intro17}.

Many hybrid models have been formulated in literature which couple the evolution of individual self--interacting endothelial cells with continuous models for chemical species diffusion and extracellular matrix evolution, highlighting the key role of chemotactic, haptotatic and contact guidance cues for the vascular network formation \cite{intro18,intro19}, together with the mechanics of cell--matrix and cell--cell interactions to describe cells migration through the extracellular matrix. 

For what concerns the modeling of tumor growth coupled with angiogenesis, an hybrid diffuse interface model of multispecies tumor growth coupled with tumor--induced angiogenesis was recently presented in \cite{intro20}, where a mixture dynamics guided by cell--cell and cell--matrix adhesion of different tumor cells species, i.e. multiple viable cells and necrotic cells, is coupled through growth and chemotaxis with a diffusible cells substrate and with a discrete random--walk angiogenesis triggered by a diffusible angiogenetic factor in the neighborhood of the tumor boundary.   

All the aforementioned models are solved on simplified geometries and do not consider the patient-specific geometry, tissue microstructure and anisotropy, which influence both the tumor cells motility, the migratory capability of endothelial cells along white matter fiber tracts and their infiltrative mechanics in the tumor tissues. Inspired by the diffuse interface mixture model introduced in \cite{intro20}, we develop a mixture model composed by a viable and a necrotic cell phases, an angiogenic phase representing tumor--induced vasculature and a liquid phase, coupled with the evolution of chemical species representing oxygen and TAF. This model is a generalization of the GBM growth model in a--vascular state, composed by a biphasic mixture, previously studied by one of the authors in \cite{agosti2}, and calibrated on patient--specific data in \cite{agosti3}. The present model is derived from variational principles which typically describe the next to equilibrium dynamics of dissipative systems composed by soft tissues of interacting constituents in isothermal situations \cite{onsager1,onsager2}, taking into account through proper biological constitutive assumptions both the tumor cells adhesion mechanics and the infiltrative mechanics of tumor--induced vasculature in the tumor tissue. Here, the endothelial cells are described through a continuum model, inspired by \cite{intro16}, since the available informations on the tumor--induced vasculature coming from neuroimaging data are at the tissue scale. The model takes the form of two degenerate Cahn--Hilliard equations with chemotaxis and source terms for the tumor cells species coupled with a Keller--Segel system for the endothelial cells and the TAF, and a reaction diffusion equation for the nutrient species. The patient--specific geometric and structural informations are included through proper constitutive assumptions in the spatial dependence of some model parameters, determined by quantitative informations extracted from neuroimaging data, while the initial distributions of tumor cells and tumor--induced vasculature are obtained by segmentation maps and post--processing maps respectively, generated from the analysis of neuroimaging data. We also include in the theory the therapy effects associated to the standard Stupp protocol, as previously done in \cite{agosti2}.

Given such a theory which describes the relation between the tumor, the vasculature and therapy, informed by neuroimaging data, it's possible to perform in silico evaluations and optimization of therapy on a patient--specific base, such as time schedule and dose optimization and dose painting \cite{intro21}. For this reason, as a proof--of--concept of the application of our model in this field, we formulate a finite element approximation of the model, which conserves the qualitative properties of the solutions of the continuous model (i.e. positivity of the concentration values and saturation of the mixture), and we perform numerical simulations on different test cases based on the data of a patient affected by GBM.

The novelties of the present work with respect to previous studies in literature are the following: 
\begin{itemize}
\item Derivation of a multiphasic tumor growth model with angiogenesis which takes into account the adhesion properties between viable and necrotic tumor cells and the extracellular matrix and the infiltrative mechanics of tumor--induced vasculature in the tumor tissues, informed by neuroimaging data giving informations about patient's brain geometry, tissue microstructure, white matter fiber orientation and tumor--induced vasculature distribution;
\item Use of state of the art processing steps to extract quantitative informations for fiber orientations from DTI data and vasculature distribution from DSCE--MRI perfusion data;
\item Coupling of the model inputs--outputs with state of the art automatic segmentation tools, based on Deep Learning, which give the distribution of different tumor components from structural MRI at different temporal stages;
\item Formulation of a finite element approximation of the model which deals with the highly non linear and degenerate nature of the system equations;
\item Study of numerical simulations for different test cases calibrated on patient--specific data, which, by comparing numerical predictions with data at later times, investigate the feasibility of the designed computational platform to possibly assist in clinical decisions.
\end{itemize}
In particular, quantitative maps for the relative Cerebral Blood Volumes (rCBV), obtained from the analysis of DSCE--MRI perfusion data taking into account the presence of blood--brain--barrier leakage induced by tumor growth, are used to reconstruct the tumor--induced vasculature distribution, due to their correlation with vasculature density \cite{intro22}. Moreover, a state of the art Deep Learning automatic segmentation tool \cite{nnunet} is used to solve the segmentation task of tumor tissues, starting from a collection of preprocessed structural MRI data, into an enhancing component, a necrotic component and edema. 

The paper is organized as follows. In Section $2$ we derive the model equations, stating the main constitutive and biological assumptions. In Section $3$ we describe the main data processing steps and we derive the finite element numerical approximation of the model. Finally in Section $4$ we deal with the model parameters selection and show numerical results, based on patient--specific data. We investigate two test cases before surgery, associated with high and low nutrient supply inside the tumor, and a test case after surgery, with the application of a standard therapy schedule, showing qualitative comparisons between numerical predictions and data. We conclude with final considerations, model limitations and future developments.

\section{Theory}
\label{sec:theory}
In this Section we derive a continuous model for tumor growth with angiogenesis, based on  a mixture model composed by a viable and a necrotic cell phases,
an angiogenic phase representing tumor--induced vasculature and a liquid phase, coupled with the evolution of chemical species representing oxygen and TAF.
The model is derived from the mass balance equations for the mixture components and the transport equations for the chemicals, given biological constitutive assumptions which comply with the second laws of thermodynamics in isothermal situations and with the least dissipation principle for next to equilibrium dynamics (which in turn imply proper momentum balance equations for the mixture). The model takes into account the chemotaxis and haptotaxis phenomena, the tumor cells adhesion mechanics and the infiltrative mechanics of tumor--induced vasculature in the tumor tissue. 

\subsection{Derivation of the model}
We consider a saturated, closed and incompressible mixture in an open bounded domain $\Omega\subset \mathbb{R}^3$, composed by a viable tumor phase $v$ with volume fraction $\phi_v$, a necrotic tumor phase $d$ with volume fraction $\phi_d$,  a liquid phase $l$ composed by liquid, healthy cells and normal vasculature, with volume fraction $\phi_l$, and an angiogenetic phase $a$ composed by tumor--induced new vasculature with volume fraction $\phi_a$. We take as an approximation that all the phases have a constant density $\gamma$, equal to the water density (since the cells are mostly composed by water). We assume that the mixture dynamics is coupled with the evolution of massless chemicals, comprising a nutrient species, with concentration (number of moles) $n$, and an angiogenetic factor, with concentration (number of moles) $c$. Each mixture component satisfies a mass continuity equation, while the massless nutrient and chemical species satisfy generic transport equations:
\begin{align}
\label{eqn:1a}
\frac{\partial \phi_v}{\partial t}+\text{div}(\phi_v \mathbf{v})+\text{div}(\mathbf{J_v})=\frac{\Gamma_v(\phi_v,\phi_d,\phi_l,a,n,c)}{\gamma},\\
\label{eqn:1b}
\frac{\partial \phi_d}{\partial t}+\text{div}(\phi_d \mathbf{v})+\text{div}(\mathbf{J_d})=\frac{\Gamma_d(\phi_v,\phi_d,\phi_l,a,n,c)}{\gamma},\\
\label{eqn:1c}
\frac{\partial \phi_a}{\partial t}+\text{div}(\phi_a \mathbf{v})+\text{div}(\mathbf{J_a})=\frac{\Gamma_a(\phi_v,\phi_d,\phi_l,a,n,c)}{\gamma},\\
\label{eqn:1d}
\frac{\partial \phi_l}{\partial t}+\text{div}(\phi_l \mathbf{v})+\text{div}(\mathbf{J_l})=\frac{\Gamma_l(\phi_v,\phi_d,\phi_l,a,n,c)}{\gamma},\\
\label{eqn:1e}
\frac{\partial n}{\partial t}+\text{div}(n \mathbf{v})+F_n=S_n(\phi_v,\phi_d,\phi_l,a,n,c),\\
\label{eqn:1f}
\frac{\partial c}{\partial t}+\text{div}(c \mathbf{v})+F_c=S_c(\phi_v,\phi_d,\phi_l,a,n,c),
\end{align}
with $\phi_v+\phi_d+\phi_a+\phi_l=1$, $\Gamma_v+\Gamma_d+\Gamma_a+\Gamma_l=0$, $\mathbf{J_v}=\phi_v(\mathbf{v_v}-\mathbf{v})$, $\mathbf{J_d}=\phi_d(\mathbf{v_d}-\mathbf{v})$, $\mathbf{J_a}=\phi_a(\mathbf{v_a}-\mathbf{v})$, $\mathbf{J_l}=\phi_l(\mathbf{v_l}-\mathbf{v})$. Here, $\mathbf{v}=\phi_v \mathbf{v_v}+\phi_d \mathbf{v_d}+\phi_a \mathbf{v_a}+\phi_l \mathbf{v_l}$ is the volume-averaged mixture velocity, which satisfies the incompressibility condition 
\begin{equation}
\label{eqn:2}
\text{div}\mathbf{v}=0,
\end{equation}
as a consequence of the saturation and the closedness properties of the mixture. 
The chemical species $n$ and $c$ are assumed to be advected by the mixture velocity $\mathbf{v}$. The terms $-F_n$ and $-F_c$ are generic transport terms to be determined in relation with the specific free energy of the system, while the source terms $S_n$ and $S_c$ represent source and consumption terms for the chemicals and must be constitutively assigned. The source terms $\Gamma_v, \Gamma_d, \Gamma_a$ represent cells proliferation, death and conversion of mass between the phases, while we take $\Gamma_l=-\left(\Gamma_v+\Gamma_d+\Gamma_a\right)$.

We make the following biological and modeling assumptions:
\begin{itemize}
\item[\textbf{A1}] The tumor cells interaction is predominant with respect to the adhesion between the tumor cells and the host tissue, and the viable and the dead tumor cells have the same adhesion properties; 
\item[\textbf{A2}] The interfacial energy between the tumor cells and the host tissues is expressed by a diffuse--interface term;
\item[\textbf{A3}] The endothelial cells of the tumor--induced vasculature constitute a self--interacting phase in the mixture (the new vasculature tips and sprouts evolve by the coupling with the angiogenetic chemical factor and by self interaction, i.e. by making loops and linings interacting with particles from the same phase);
\item[\textbf{A4}] The viable tumor cells can migrate to regions with higher nutrient concentration, being coupled to the nutrient by a chemotactic term;
\item[\textbf{A5}] The endothelial cells of the tumor--induced vasculature can migrate to regions with higher angiogenetic factor concentration, being coupled to the angiogenetic factor by a chemotactic term;
\item[\textbf{A6}] The chemical species are characterized by random molecular motion throughout the mixture; 
\item[\textbf{A7}] The predominant friction forces in the relative motion between the mixture constituents are the friction between the tumor cells and the liquid phase on one hand and the friction between the endothelial cells and both the liquid and the tumor phases on the other hand.
\end{itemize}
With assumptions \textbf{A1}--\textbf{A6}, the free energy of the system, given by its internal energy minus its entropy, takes the following general form

\begin{align}
\label{eqn:3}
&E(\phi_v,\phi_d,a,n,c)=\int_{\Omega}e(\phi_v,\phi_d,a,n,c)\,d\mathbf{x}=\\
\notag &\int_{\Omega}\Pi \biggl(k_T\phi_T\log\phi_T+k_T(1-\phi_T)\log(1-\phi_T)+k_a\phi_a(\log(\phi_a)-1)+k_I\phi_T(1-\phi_T)+\\
& \notag \frac{\epsilon^2}{2}|\nabla \phi_T|^2+\frac{\mathbf{D_n}}{2}\nabla n \cdot \nabla n-\chi_vn\phi_v+\frac{\mathbf{D_c}}{2}\nabla c \cdot \nabla c-\chi_a\phi_ac\biggr)d\mathbf{x},
\end{align}
where $e(\phi_v,\phi_d,a,n,c)$ is the free energy per unit volume, and $\phi_T=\phi_v+\phi_d$ is the total tumor cell concentration. $\Pi$ is the Young modulus of the brain tissue, in units of $[Pa]$, while $k_T, k_a, k_I$ are specific coefficients for each component's contribution. The term $-k_T\phi_T\log\phi_T-k_T(1-\phi_T)\log(1-\phi_T)$ is the entropy associated to the binary interaction between the tumor cells and their surrounding, $k_I\phi_T(1-\phi_T)$ is the internal energy associated to the interaction between  the tumor cells and their surrounding, and  $-k_a\phi_a(\log(\phi_a)-1)$ is the entropy associated to the self--interacting cells of the phase $a$ (see e.g. \cite{suzuki} for the expression of the entropy for a system of self--interacting particles). The term $\frac{\epsilon^2}{2}|\nabla \phi_T|^2$ represents the interfacial internal energy between the tumor cells and the host tissue, with $\epsilon$, in units of $[m]$, the interfacial thickness. The terms $\frac{\mathbf{D_n}}{2}\nabla n \cdot \nabla n$ and $\frac{\mathbf{D_c}}{2}\nabla c \cdot \nabla c$ represent the contribution to the internal energy from the random motion of the chemical species resulting in diffusive behaviors along concentration gradients, where $\mathbf{D_n}$ and $\mathbf{D_c}$ are the matrices of mean deviation of the displacement of the particles in the different directions, with units of $[mm^2/\text{Mol}^2]$, possibly representing anisotropic diffusion in the tissues. Finally, $-\chi_vn\phi_v$ and $-\chi_a\phi_ac$ are interaction terms associated to chemotaxis, with chemotactic coefficients $\chi_v$ and $\chi_a$ in units of $[\text{Mol}^{-1}]$. According to assumption \textbf{A1} and following e.g. \cite{agosti1,agosti2}, we substitute the expression of the homogeneous free energy associated to the binary interaction between the tumor cells and the surrounding tissues by a phenomenological term for the cell-cell mechanical interactions, which has the form of a single well potential of the \textit{Lennard--Jones} type:
\begin{equation}
\label{eqn:4}
\psi(\phi_T)=-(1-\bar{\phi})\log(1-\phi_T)-\frac{\phi_T^3}{3}-(1-\bar{\phi})\left(\frac{\phi_T^2}{2}+\phi_T\right),
\end{equation}
where $\bar{\phi}$ represents the equilibrium value of the tumor cell concentration at which no interacting force is exerted between the cells. Hence, the free energy of the system takes the form 
\begin{align}
\label{eqn:5}
&E(\phi_v,\phi_d,\phi_l,a,n,c)=\\
\notag &\int_{\Omega}\Pi \biggl(\psi(\phi_T)+k_a\phi_a(\log(\phi_a)-1)+\frac{\epsilon^2}{2}|\nabla \phi_T|^2+\frac{\mathbf{D_n}}{2}\nabla n \cdot \nabla n-\chi_vn\phi_v+\frac{\mathbf{D_c}}{2}\nabla c \cdot \nabla c-\chi_a\phi_ac\biggr)d\mathbf{x}.
\end{align}

We now derive the model for the dynamics of the system satisfying the second law of thermodynamics in isothermal situations and a generalized variational principle of least dissipation.\\
Let us firstly introduce the functional derivatives of the free energy $E(\phi_v,\phi_d,a,n,c)$ (we take $\Pi=k_a=1$ in the model derivation for simplicity):
\begin{align*} 
&\mu_v:=\frac{\delta E}{\delta \phi_v}=\psi'(\phi_T)-\epsilon^2 \Delta \phi_T-\chi_{v}n,\\
&\mu_d:=\frac{\delta E}{\delta \phi_d}=\psi'(\phi_T)-\epsilon^2 \Delta \phi_T,\\
&\mu_a:=\frac{\delta E}{\delta \phi_a}=\log(\phi_a)-\chi_{a}c,\\
&\eta:=\frac{\delta E}{\delta n}=-\text{div}(\mathbf{D_n}\nabla n)-\chi_v\phi_v,\\
&\theta:=\frac{\delta E}{\delta c}=-\text{div}(\mathbf{D_c}\nabla c)-\chi_a\phi_a.
\end{align*} 
We multiply equation \eqref{eqn:1a} by $\mu_v$, equation \eqref{eqn:1b} by $\mu_d$, equation \eqref{eqn:1c} by $\mu_a$, equation \eqref{eqn:1e} by $\eta$, equation \eqref{eqn:1f} by $\theta$, \eqref{eqn:2} by a multiplier $p$, integrate the resulting terms over a generic volume $R(t)\subset \Omega$ transported with the mixture velocity and sum the three contributions. We get
\begin{align*}
& \int_{R(t)}\mu_v \left(\frac{\partial \phi_v}{\partial t}+\mathbf{v}\cdot \nabla \phi_v + \phi_v \text{div}\mathbf{v}+\text{div}\mathbf{J}_v\right)\,d\mathbf{x}+\\
& \int_{R(t)}\mu_d \left(\frac{\partial \phi_d}{\partial t}+\mathbf{v}\cdot \nabla \phi_d + \phi_d \text{div}\mathbf{v}+\text{div}\mathbf{J}_d\right)\,d\mathbf{x}+\\
& \int_{R(t)}\mu_a \left(\frac{\partial \phi_a}{\partial t}+\mathbf{v}\cdot \nabla \phi_a + \phi_a \text{div}\mathbf{v}+\text{div}\mathbf{J}_a\right)\,d\mathbf{x}+\\
& \int_{R(t)}\eta \left(\frac{\partial n}{\partial t}+\mathbf{v}\cdot \nabla n + n \text{div}\mathbf{v}+F_n\right)\,d\mathbf{x}+\\
&\int_{R(t)}\theta \left(\frac{\partial c}{\partial t}+\mathbf{v}\cdot \nabla c + c \text{div}\mathbf{v}+F_c\right)\,d\mathbf{x}+\int_{R(t)}p\text{div}\mathbf{v}\,d\mathbf{x}=\\
& \int_{R(t)}\left(\frac{\Gamma_v}{\gamma}\mu_v+\frac{\Gamma_d}{\gamma}\mu_d+\frac{\Gamma_a}{\gamma}\mu_a+S_n\eta+S_c\theta\right)\,d\mathbf{x}.
\end{align*}
Inserting the explicit expressions of $\mu_v$, $\mu_d$, $\mu_a$, $\eta$ and $\theta$, rearranging the terms and integrating by parts we get
\begin{align*}
&\int_{R(t)}\biggl\{\left(\psi'(\phi_T)-\epsilon^2 \Delta\phi_T\right)\left(\frac{\partial \phi_T}{\partial t}+\mathbf{v}\cdot \nabla \phi_T\right)-\chi_vn\left(\frac{\partial \phi_v}{\partial t}+\mathbf{v}\cdot \nabla \phi_v\right)-\\
& \left(\text{div}(\mathbf{D_n}\nabla n)+\chi_v\phi_v\right)\left(\frac{\partial n}{\partial t}+\mathbf{v}\cdot \nabla n\right)+\left(\log(\phi_a)-\chi_{a}c\right)\left(\frac{\partial \phi_a}{\partial t}+\mathbf{v}\cdot \nabla \phi_a\right)-\\
&\left(\text{div}(\mathbf{D_c}\nabla c)+\chi_a\phi_a\right)\left(\frac{\partial c}{\partial t}+\mathbf{v}\cdot \nabla c\right)+\left(\mu_v\phi_v+\mu_d\phi_d+\mu_a\phi_a+\eta n+\theta c+p\right)\text{div}\mathbf{v}-\\
&\left(\mathbf{J}_v\cdot \nabla \mu_v+\mathbf{J}_d\cdot \nabla \mu_d+\mathbf{J}_a\cdot \nabla \mu_a-F_n\eta -F_c\theta \right)\biggr\}\,d\mathbf{x}+\int_{\partial R(t)}\left(\mu_v\mathbf{J}_v+\mu_d\mathbf{J}_d+\mu_a\mathbf{J}_a\right)\cdot \boldsymbol{\nu}\,dS=\\
&\int_{R(t)}\left(\frac{\Gamma_v}{\gamma}\mu_v+\frac{\Gamma_d}{\gamma}\mu_d+\frac{\Gamma_a}{\gamma}\mu_a+S_n\eta+S_c\theta\right)\,d\mathbf{x},
\end{align*}
where $\boldsymbol{\nu}$ is the outer normal to $\partial R(t)$.
Using the identity
\[
0=\int_R\left(\mathbf{M}\nabla f \cdot \left[\nabla \nabla f\right] \mathbf{v}+\frac{\left[\nabla \mathbf{M}\right]^T}{2}\nabla f \otimes \nabla f \cdot \mathbf{v}-\nabla\left(\frac{\mathbf{M}}{2}\nabla f \cdot \nabla f\right)\cdot \mathbf{v}\right)\,d\mathbf{x},
\]
which is valid for any sufficiently regular matrix $\mathbf{M}$ and scalar function $f$, rearranging terms and integrating by parts, we obtain
\begin{align*}
&\int_{R(t)}\biggl\{\psi'(\phi_T)\left(\frac{\partial \phi_T}{\partial t}+\mathbf{v}\cdot \nabla \phi_T\right)+\epsilon^2\nabla \phi_T\cdot \left(\frac{\partial \nabla \phi_T}{\partial t}+\left[\nabla \nabla \phi_T\right]\mathbf{v}\right)-\left((\epsilon^2 \Delta\phi_T-\psi'(\phi_T))\nabla \phi_T\cdot \mathbf{v}\right)\\
&-\left(\chi_nn\nabla \phi_v\cdot \mathbf{v}\right)-\left(\nabla\left(\epsilon^2\frac{|\nabla \phi_T|^2}{2}+\psi(\phi_T)\right)\cdot \mathbf{v}\right)-\chi_vn\left(\frac{\partial \phi_v}{\partial t}+\mathbf{v}\cdot \nabla \phi_v\right)+\\
&\mathbf{D_n}\nabla n\cdot \left(\frac{\partial \nabla n}{\partial t}+\left[\nabla \nabla n\right]\mathbf{v}\right)-\left(\left(\text{div}(\mathbf{D_n}\nabla n)+\chi_v\phi_v\right)\nabla n\cdot \mathbf{v}\right)-\left(\nabla\left(\frac{\mathbf{D}_n}{2}\nabla n \cdot \nabla n-\chi_nn\phi_v\right)\cdot \mathbf{v}\right)\\
&+\left(\frac{\left[\nabla \mathbf{D}_n\right]^T}{2}\nabla n \otimes \nabla n \cdot \mathbf{v}\right)-\chi_v\phi_v\left(\frac{\partial n}{\partial t}+\mathbf{v}\cdot \nabla n\right)+\left(\log(\phi_a)-\chi_{a}c\right)\left(\frac{\partial \phi_a}{\partial t}+\mathbf{v}\cdot \nabla \phi_a\right)\\
&+\mathbf{D_c}\nabla c\cdot \left(\frac{\partial \nabla c}{\partial t}+\left[\nabla \nabla c\right]\mathbf{v}\right)-\left(\left(\text{div}(\mathbf{D_c}\nabla c)+\chi_a\phi_a\right)\nabla c\cdot \mathbf{v}\right)-\left(\nabla\left(\frac{\mathbf{D}_c}{2}\nabla c \cdot \nabla c-\chi_ac\phi_a\right)\cdot \mathbf{v}\right)\\
&+\left(\frac{\left[\nabla \mathbf{D}_c\right]^T}{2}\nabla c \otimes \nabla c \cdot \mathbf{v}\right)+\left(\log(\phi_a)-\chi_{a}c\right)\nabla \phi_a\cdot \mathbf{v}-\left(\nabla\left(\phi_a(\log(\phi_a)-1)\right)\cdot \mathbf{v}\right)\\
&-\chi_a\phi_a\left(\frac{\partial c}{\partial t}+\mathbf{v}\cdot \nabla c\right)+\left(\mu_v\phi_v+\mu_d\phi_d+\mu_a\phi_a+\eta n+\theta c+p\right)\text{div}\mathbf{v}-\\
&\left(\mathbf{J}_v\cdot \nabla \mu_v+\mathbf{J}_d\cdot \nabla \mu_d+\mathbf{J}_a\cdot \nabla \mu_a-F_n\eta -F_c\theta \right)\biggr\}\,d\mathbf{x}+\int_{\partial R(t)}\left(\mu_v\mathbf{J}_v+\mu_d\mathbf{J}_d+\mu_a\mathbf{J}_a\right)\cdot \boldsymbol{\nu}\,dS-\\
&\int_{\partial R(t)}\left(\epsilon^2\left(\frac{\partial \phi_T}{\partial t}\right)\nabla \phi_T+\left(\frac{\partial n}{\partial t}\right)\mathbf{D_n} \nabla n+\left(\frac{\partial c}{\partial t}\right)\mathbf{D_c} \nabla c\right)\cdot \boldsymbol{\nu}\,dS=\\
& \int_{R(t)}\left(\frac{\Gamma_v}{\gamma}\mu_v+\frac{\Gamma_d}{\gamma}\mu_d+\frac{\Gamma_a}{\gamma}\mu_a+S_n\eta+S_c\theta\right)\,d\mathbf{x}.
\end{align*}
Recalling the form of the free energy density $e$ in \eqref{eqn:5}, we finally obtain
\begin{align}
\label{eqn:8}
&\int_{R(t)}\biggl(\frac{\partial e}{\partial t}+\mathbf{v}\cdot \nabla e + e\text{div}\mathbf{v}+\mathbf{v}\cdot\left(\mu_v \nabla \phi_v+\mu_d \nabla \phi_d+\mu_a \nabla \phi_a+\eta \nabla n + \theta \nabla c\right)-\\
&\notag \mathbf{v}\cdot\left(\nabla \left(p+\mu_v\phi_v+\mu_d\phi_d+\mu_a\phi_a+\eta n+\theta c\right)\right)-\\
&\notag\left(\mathbf{J}_v\cdot \nabla \mu_v+\mathbf{J}_d\cdot \nabla \mu_d+\mathbf{J}_a\cdot \nabla \mu_a-F_n\eta -F_c\theta \right)\biggr)\,d\mathbf{x}+\\
&\notag \int_{\partial R(t)}\left(\mu_v\mathbf{J}_v+\mu_d\mathbf{J}_d+\mu_a\mathbf{J}_a\right)\cdot \boldsymbol{\nu}\,dS-\\
&\notag \int_{\partial R(t)}\left(\epsilon^2\left(\frac{\partial \phi_T}{\partial t}\right)\nabla \phi_T+\left(\frac{\partial n}{\partial t}\right)\mathbf{D_n} \nabla n+\left(\frac{\partial c}{\partial t}\right)\mathbf{D_c} \nabla c\right)\cdot \boldsymbol{\nu}\,dS+\\
&\notag \int_{\partial R(t)}\left(\mu_v\phi_v+\mu_d\phi_d+\mu_a\phi_a+\eta n+\theta c+p-e\right)\mathbf{v}\cdot \boldsymbol{\nu}\,dS=\\
&\notag\int_{R(t)}\left(\frac{\Gamma_v}{\gamma}\mu_v+\frac{\Gamma_d}{\gamma}\mu_d+\frac{\Gamma_a}{\gamma}\mu_a+S_n\eta+S_c\theta\right)\,d\mathbf{x}.
\end{align}
The second law of thermodynamics in isothermal situations and with source terms takes the form of the following dissipation inequality \cite{garcke1,gurtin1}
\begin{equation}
\label{eqn:9}
\frac{d}{dt}\int_{R(t)}e\leq -\int_{\partial R(t)}\mathbf{J_E}\cdot \boldsymbol{\nu}dS+\int_{R(t)}\biggl(\frac{\Gamma_v}{\gamma}m_v+\frac{\Gamma_d}{\gamma}m_d+\frac{\Gamma_a}{\gamma}m_a+S_nm_n+S_cm_c\biggr)d\mathbf{x},
\end{equation}
for each material volume $R(t)\subset \Omega$, with the energy flux $\mathbf{J_E}$ and the multipliers $m_v, m_d, m_a, m_n, m_c$ to be determined.
By comparing \eqref{eqn:8} with \eqref{eqn:9}, the following constitutive assumptions can be made in order for the system to fulfill the dissipation inequality \eqref{eqn:9}
\begin{align}
\label{eqn:10a}
\bar{p}=&p+\mu_v\phi_v+\mu_d\phi_d+\mu_a\phi_a+\eta n+\theta c,\\
\label{eqn:10b}
\mathbf{v}=&-k\left(\nabla \bar{p}-\left(\mu_v \nabla \phi_v+\mu_d \nabla \phi_d+\mu_a \nabla \phi_a+\eta \nabla n + \theta \nabla c\right)\right),\\
\label{eqn:10c}
\mathbf{J_v}=&-b_v(\phi_v,\phi_d,\phi_a)\mathbf{T}\nabla \mu_v,\\
\label{eqn:10d}
\mathbf{J_d}=&-b_d(\phi_v,\phi_d,\phi_a)\mathbf{T}\nabla \mu_d,\\
\label{eqn:10e}
\mathbf{J_a}=&-b_a(\phi_v,\phi_d,\phi_a)\mathbf{T}\nabla \mu_a,\\
\label{eqn:10f}
F_n=&\alpha_n\eta,\\
\label{eqn:10g}
F_c=&\alpha_c\theta,\\
\label{eqn:10h}
\mathbf{J_E}=&\mu_v\mathbf{J}_v+\mu_d\mathbf{J}_d+\mu_a\mathbf{J}_a+\epsilon^2\left(\frac{\partial \phi_T}{\partial t}\right)\nabla \phi_T+\left(\frac{\partial n}{\partial t}\right)\mathbf{D_n} \nabla n+\left(\frac{\partial c}{\partial t}\right)\mathbf{D_c} \nabla c\\
& \notag +\left(\mu_v\phi_v+\mu_d\phi_d+\mu_a\phi_a+\eta n+\theta c+p-e\right)\mathbf{v},\\
m_v=&\mu_v, m_d=\mu_d, m_a=\mu_a, m_n=\eta, m_c=\theta,
\end{align}
where $k$ is a positive friction parameter, $\alpha_n,\alpha_c$ are positive coefficients related to the time scales of the dynamics of the chemical species, in units of $[\text{Mol}^2/ Pa \, \text{day}]$, $b_v, b_d, b_a$ are positive mobilities, assumed to be dependent only on the mixture components concentrations, and $\mathbf{T}$ is a symmetric positive definite tensor, representing anisotropy in the cells motility.
Inserting \eqref{eqn:10a}--\eqref{eqn:10g} in \eqref{eqn:1a}--\eqref{eqn:1f}, we get the following system of equations
\begin{equation}
\label{eqn:11}
\begin{cases}
\mathbf{v}=-k\left(\nabla \bar{p}-\left(\mu_v \nabla \phi_v+\mu_d \nabla \phi_d+\mu_a \nabla \phi_a+\eta \nabla n + \theta \nabla c\right)\right),\\
\text{div}\mathbf{v}=0,\\
\frac{\partial \phi_v}{\partial t}+\mathbf{v}\cdot \nabla \phi_v - \text{div}\biggl(b(\phi_v)\mathbf{T}\nabla \mu_v\biggr)=\frac{\Gamma_v}{\gamma} ,\\
\frac{\partial \phi_d}{\partial t}+\mathbf{v}\cdot \nabla \phi_d - \text{div}\biggl(b(\phi_d)\mathbf{T}\nabla \mu_d\biggr)=\frac{\Gamma_d}{\gamma} ,\\
\frac{\partial \phi_a}{\partial t}+\mathbf{v}\cdot \nabla \phi_a - \text{div}\biggl(b(\phi_a)\mathbf{T}\nabla \mu_a\biggr)=\frac{\Gamma_a}{\gamma} ,\\
\mu_v=\psi'(\phi_T)-\epsilon^2 \Delta \phi_T-\chi_v n,\\
\mu_d=\psi'(\phi_T)-\epsilon^2 \Delta \phi_T,\\
\mu_a=\log(\phi_a)-\chi_{a}c,\\
\frac{\partial n}{\partial t}+\mathbf{v}\cdot \nabla n - \alpha_n\text{div}\left(\mathbf{D}_n\nabla n\right)-\alpha_n\chi_v \phi_v= S_n,\\
\frac{\partial c}{\partial t}+\mathbf{v}\cdot \nabla c - \alpha_c\text{div}\left(\mathbf{D}_c\nabla c\right)-\alpha_c\chi_a \phi_a= S_c,
\end{cases}
\end{equation} 
in $\Omega$, which we endow with the homogeneous boundary conditions
\begin{align}
\label{eqn:11bis}
&b_v\mathbf{T}\nabla \mu_v \cdot \boldsymbol{\nu}|_{\partial \Omega}=b_d\mathbf{T}\nabla \mu_d \cdot \boldsymbol{\nu}|_{\partial \Omega}=b_a\mathbf{T}\nabla \mu_a \cdot \boldsymbol{\nu}|_{\partial \Omega}=\nabla \phi_T \cdot \boldsymbol{\nu}|_{\partial \Omega}=\nabla n \cdot \boldsymbol{\nu}|_{\partial \Omega}=\\
& \notag\nabla c \cdot \boldsymbol{\nu}|_{\partial \Omega}=\mathbf{v}|_{\partial \Omega} =0,
\end{align}
which imply that
\[
\mathbf{J_E}\cdot \boldsymbol{\nu}|_{\partial \Omega}=0.
\]
A solution of system \eqref{eqn:11}, supplemented with the boundary conditions \eqref{eqn:11bis}, formally satisfies the following energy equality
\begin{align}
\label{eqn:12}
&\frac{dE}{dt}+\frac{1}{k}\int_{\Omega}\mathbf|v|^2 d\mathbf{x} +\int_{\Omega}\left(b_v\mathbf{T}\nabla \mu_v\cdot \nabla \mu_v+b_d\mathbf{T}\nabla \mu_d\cdot \nabla \mu_d+b_a\mathbf{T}\nabla \mu_a\cdot \nabla \mu_a+\alpha_n\eta^2+\alpha_c\theta^2\right) d\mathbf{x}\\
& \notag =\int_{\Omega}\left(\frac{\Gamma_v}{\gamma}\mu_v+\frac{\Gamma_d}{\gamma}\mu_d+\frac{\Gamma_a}{\gamma}\mu_a+S_n\eta+S_c\theta\right)\,d\mathbf{x},
\end{align}
which is obtained by multiplying the first equation of \eqref{eqn:11} by $\frac{1}{k}\mathbf{v}$, the second equation by $\bar{p}$, the third, fourth and fifth equations by $\mu_v, \mu_d, \mu_a$ respectively, the sixth and seventh equations by $\eta, \theta$ respectively, integrating over $\Omega$ and summing all the contributions.

In order to close the system \eqref{eqn:11} and determine particular forms for the mobility functions $b_v,b_d,b_a$, we apply the Onsager Variational Principle (OVP) \cite{onsager1}, which defines the irreversible non--equilibrium dynamics for near--equilibrium thermodynamic systems in terms of linear fluxes--forces balance equations. This principle has been widely applied for the derivation of continuum phenomenological models of soft matter, see e. g. \cite{onsager2,onsager3,onsager4}, based on the observation that the macroscopic behavior of soft matter is driven by the interactions between its consituents units, at the constituent scale, which drive the system locally in a metastable regime out of equilibrium, and by the dissipation mechanisms operating in the system itself. In isothermal situations, the OVP takes the following form: given a set of slow state variables $\mathbf{x}_i, i=1, \dots, n$, the dynamics of the system is described by the thermodynamic fluxes which minimize the Onsager functional $\mathcal{O}(\dot{\mathbf{x}}_i)=\Phi(\dot{\mathbf{x}}_i)+\dot{E}(\mathbf{x}_i,\dot{\mathbf{x}}_i)$, where $\Phi$ is the dissipation function, which is quadratic in $\dot{\mathbf{x}}_i$ as a near--equilibrium approximation, and $E$ is the free energy of the system. In our case, we minimize \eqref{eqn:12} with respect to the variables $\mathbf{v}_v, \mathbf{v}_d, \mathbf{v}_a, \mathbf{v}_l$, thus obtaining the momentum balance equations for the four phases of the mixture as linear fluxes--forces relations, to be supplemented to the mass balance equations \eqref{eqn:1a}-\eqref{eqn:1d} in the mixture dynamics description. We thus rewrite \eqref{eqn:12} as 
\begin{align}
\label{eqn:13}
&\int_{\Omega}\frac{\delta E}{\delta \phi_v}\biggl(\frac{\Gamma_v}{\gamma}-\text{div}(\phi_v \mathbf{v}_v)\biggr)d\mathbf{x}+\int_{\Omega}\frac{\delta E}{\delta \phi_d}\biggl(\frac{\Gamma_d}{\gamma}-\text{div}(\phi_d \mathbf{v}_d)\biggr)d\mathbf{x}+\int_{\Omega}\frac{\delta E}{\delta \phi_a}\biggl(\frac{\Gamma_a}{\gamma}-\text{div}(\phi_a \mathbf{v}_a)\biggr)d\mathbf{x}\\
& \notag +\int_{\Omega}\frac{\delta E}{\delta n}\biggl(S_n-\text{div}(n \mathbf{v})-F_n\biggr)d\mathbf{x}+\int_{\Omega}\frac{\delta E}{\delta c}\biggl(S_c-\text{div}(c \mathbf{v})-F_c\biggr)d\mathbf{x}+\frac{1}{k}\int_{\Omega}\mathbf|v|^2 d\mathbf{x}\\
& \notag+\int_{\Omega}\left(b_v\mathbf{T}\nabla \mu_v\cdot \nabla \mu_v+b_d\mathbf{T}\nabla \mu_d\cdot \nabla \mu_d+b_a\mathbf{T}\nabla \mu_a\cdot \nabla \mu_a+\alpha_n\eta^2+\alpha_c\theta^2\right) d\mathbf{x}\\
& \notag =\int_{\Omega}\left(\frac{\Gamma_v}{\gamma}\mu_v+\frac{\Gamma_d}{\gamma}\mu_d+\frac{\Gamma_a}{\gamma}\mu_a+S_n\eta+S_c\theta\right)\,d\mathbf{x}
\end{align}
which gives 
\begin{align}
\label{eqn:14}
&\int_{\Omega}\mu_v \left(-\text{div}(\phi_v \mathbf{v}_v)\right)d\mathbf{x}+\int_{\Omega}\mu_d \left(-\text{div}(\phi_d \mathbf{v}_d)\right)d\mathbf{x}+\int_{\Omega}\mu_a \left(-\text{div}(\phi_a \mathbf{v}_a)\right)d\mathbf{x}+\\
& \notag\int_{\Omega}\eta \left(-\text{div}(n \mathbf{v})\right)d\mathbf{x}+\int_{\Omega}\theta \left(-\text{div}(c \mathbf{v})\right)d\mathbf{x}+\frac{1}{k}\int_{\Omega}\mathbf|v|^2 d\mathbf{x}+\\
& \notag\int_{\Omega}\left(b_v\mathbf{T}\nabla \mu_v\cdot \nabla \mu_v+b_d\mathbf{T}\nabla \mu_d\cdot \nabla \mu_d+b_a\mathbf{T}\nabla \mu_a\cdot \nabla \mu_a\right) d\mathbf{x}=0.
\end{align}
As in the framework of the OVP, we assume that the dissipation term due to the viscous interactions between the phases can be assumed to depend only on the drag between the phases and can be written as a quadratic term in the relative velocity between the phases. Following assumption \textbf{A7}, we can write
\begin{align}
\label{eqn:15a}
&\int_{\Omega}b_v(\phi_v,\phi_d,\phi_a)\mathbf{T}\nabla \mu_v \cdot \nabla \mu_v d\mathbf{x}\equiv \int_{\Omega}\frac{\phi_v}{2}M_v\mathbf{T}^{-1}(\mathbf{v_v}-\mathbf{v_l})\cdot (\mathbf{v_v}-\mathbf{v_l}) d\mathbf{x},\\
\label{eqn:15b}
&\int_{\Omega}b_d(\phi_v,\phi_d,\phi_a)\mathbf{T}\nabla \mu_d \cdot \nabla \mu_d d\mathbf{x}\equiv \int_{\Omega}\frac{\phi_d}{2}M_d\mathbf{T}^{-1}(\mathbf{v_d}-\mathbf{v_l})\cdot (\mathbf{v_d}-\mathbf{v_l}) d\mathbf{x},\\
\label{eqn:15c}
&\int_{\Omega}b_a(\phi_v,\phi_d,\phi_a)\mathbf{T}\nabla \mu_a \cdot \nabla \mu_a d\mathbf{x}\equiv \int_{\Omega}\biggl(\frac{\phi_a}{2}M_{al}\mathbf{T}^{-1}(\mathbf{v_a}-\mathbf{v_l})\cdot (\mathbf{v_a}-\mathbf{v_l}) +\\
& \notag\frac{\phi_a}{2}M_{av}\mathbf{T}^{-1}(\mathbf{v_a}-\mathbf{v_v})\cdot (\mathbf{v_a}-\mathbf{v_v})+\frac{\phi_a}{2}M_{ad}\mathbf{T}^{-1}(\mathbf{v_a}-\mathbf{v_d})\cdot (\mathbf{v_a}-\mathbf{v_d}) \biggr)d\mathbf{x},
\end{align}
where $M_{v}, M_{d}, M_{al}, M_{av}, M_{ad}$ are different friction parameters relative to the specific filtration processes.
Sobstituting \eqref{eqn:15a}--\eqref{eqn:15c} in \eqref{eqn:14}, we find that a minimum of \eqref{eqn:14} with respect to the phase velocities $\mathbf{v_v}, \mathbf{v_d}, \mathbf{v_a}$ and $\mathbf{v_l}$ satisfies the following first-order conditions
\begin{equation}
\label{eqn:16}
\begin{cases}
\phi_v \nabla \mu_v +\phi_v n \nabla \eta +\phi_v c \nabla \theta + 2\frac{\phi_v}{k} \mathbf{v}+\frac{\phi_v}{M_v}\mathbf{T}^{-1}(\mathbf{v_v}-\mathbf{v_l})-\frac{\phi_a}{M_{av}}\mathbf{T}^{-1}(\mathbf{v_a}-\mathbf{v_v})=0,\\ \\
\phi_d \nabla \mu_d +\phi_d n \nabla \eta +\phi_d c \nabla \theta + 2\frac{\phi_d}{k} \mathbf{v}+\frac{\phi_d}{M_d}\mathbf{T}^{-1}(\mathbf{v_d}-\mathbf{v_l})-\frac{\phi_a}{M_{ad}}\mathbf{T}^{-1}(\mathbf{v_a}-\mathbf{v_d})=0,\\ \\
\phi_a \nabla \mu_a +\phi_a n \nabla \eta +\phi_a c \nabla \theta + 2\frac{\phi_a}{k} \mathbf{v}+\frac{\phi_a}{M_{al}}\mathbf{T}^{-1}(\mathbf{v_a}-\mathbf{v_l})+\frac{\phi_a}{M_{av}}\mathbf{T}^{-1}(\mathbf{v_a}-\mathbf{v_v})+\\ 
\frac{\phi_a}{M_{ad}}\mathbf{T}^{-1}(\mathbf{v_a}-\mathbf{v_d})=0,\\ \\
\phi_l n \nabla \eta +\phi_l c \nabla \theta + 2\frac{\phi_l}{k} \mathbf{v}-\frac{\phi_v}{M_v}\mathbf{T}^{-1}(\mathbf{v_v}-\mathbf{v_l})-\frac{\phi_d}{M_d}\mathbf{T}^{-1}(\mathbf{v_d}-\mathbf{v_l})-\frac{\phi_a}{M_{al}}\mathbf{T}^{-1}(\mathbf{v_a}-\mathbf{v_l})=0.
\end{cases}
\end{equation}
Summing the last two equations in \eqref{eqn:16}, we get
\begin{align}
\label{eqn:16bis}
& \frac{2}{k}\mathbf{v}=-\frac{\phi_a}{\phi_a+\phi_l}\nabla \mu_a-n\nabla \eta-c\nabla \theta-\frac{\phi_a}{M_{av}(\phi_a+\phi_l)}\mathbf{T}^{-1}(\mathbf{v_a}-\mathbf{v_v})-\\
& \notag \frac{\phi_a}{M_{ad}(\phi_a+\phi_l)}\mathbf{T}^{-1}(\mathbf{v_a}-\mathbf{v_d})+\frac{\phi_v}{M_v(\phi_a+\phi_l)}\mathbf{T}^{-1}(\mathbf{v_v}-\mathbf{v_l})+\frac{\phi_d}{M_d(\phi_a+\phi_l)}\mathbf{T}^{-1}(\mathbf{v_d}-\mathbf{v_l}).
\end{align}
Inserting \eqref{eqn:16bis} in the first three equations of \eqref{eqn:16}, and summing the contributions, we get
\begin{align}
\label{eqn:16tris}
&\phi_v\nabla \mu_v+\phi_d\nabla \mu_d+\phi_a\nabla \mu_a -\frac{\phi_a(1-\phi_l)}{\phi_a+\phi_l}\nabla \mu_a+\frac{\phi_v(1+\phi_a)}{M_v(\phi_a+\phi_l)}\mathbf{T}^{-1}(\mathbf{v_v}-\mathbf{v_l})+\\
& \notag \frac{\phi_d(1+\phi_a)}{M_d(\phi_a+\phi_l)}\mathbf{T}^{-1}(\mathbf{v_d}-\mathbf{v_l})+\frac{\phi_a}{M_{al}}\mathbf{T}^{-1}(\mathbf{v_a}-\mathbf{v_l})-\frac{\phi_a(1-\phi_l)}{M_{av}(\phi_a+\phi_l)}\mathbf{T}^{-1}(\mathbf{v_a}-\mathbf{v_v})-\\
& \notag \frac{\phi_a(1-\phi_l)}{M_{ad}(\phi_a+\phi_l)}\mathbf{T}^{-1}(\mathbf{v_a}-\mathbf{v_d})=0
\end{align}
A solution of \eqref{eqn:16tris}, compatible with \eqref{eqn:15a}--\eqref{eqn:15c}, is
\begin{equation}
\label{eqn:17}
\begin{cases}
\mathbf{v_v}-\mathbf{v_l}=-\frac{1-\phi_v-\phi_d}{M_v(1+\phi_a)}\mathbf{T}\nabla \mu_v,\\
\mathbf{v_d}-\mathbf{v_l}=-\frac{1-\phi_v-\phi_d}{M_d(1+\phi_a)}\mathbf{T}\nabla \mu_d,\\
\mathbf{v_a}-\mathbf{v_l}=-\frac{1}{M_{al}}\mathbf{T}\nabla \mu_a,\\
\mathbf{v_a}-\mathbf{v_v}=-\frac{1}{M_{av}}\mathbf{T}\nabla \mu_a,\\
\mathbf{v_a}-\mathbf{v_d}=-\frac{1}{M_{ad}}\mathbf{T}\nabla \mu_a.
\end{cases}
\end{equation}
Comparing \eqref{eqn:17} with \eqref{eqn:15a}--\eqref{eqn:15c}, we get that 
\begin{align}
\label{eqn:18}
&b_v(\phi_v,\phi_d,\phi_a)=\frac{\phi_v(1-\phi_v-\phi_d)^2}{2M_v(1+\phi_a)^2}, \; b_d(\phi_v,\phi_d,\phi_a)=\frac{\phi_d(1-\phi_v-\phi_d)^2}{2M_d(1+\phi_a)^2},\\ & \notag b_a(\phi_v,\phi_d,\phi_a)=\frac{\phi_a}{2}\left(\frac{1}{M_{al}+M_{av}+M_{ad}}\right).
\end{align}
In the limit of high viscosity of the mixture, which corresponds to a low value of the friction parameter $k$, (which is appropriate for the description of the tumor dynamics, see e.g. \cite{agosti2}), we can simplify system \eqref{eqn:11} by taking $\mathbf{v}\equiv 0$. Sobstituting \eqref{eqn:18} in \eqref{eqn:11}, and indicating $M_a:=M_{al}+M_{av}+M_{ad}$ , we finally obtain the system
\begin{equation}
\label{eqn:19}
\begin{cases}
\displaystyle
\frac{\partial \phi_v}{\partial t}- \text{div}\left(\frac{\phi_v(1-\phi_v-\phi_d)^2}{2M_v(1+\phi_a)^2}\mathbf{T}\nabla \mu \right)+\text{div}\left(\chi_v\frac{\phi_v(1-\phi_v-\phi_d)^2}{2M_v(1+\phi_a)^2}\mathbf{T}\nabla n\right)=\frac{\Gamma_v}{\gamma} ,\\
\displaystyle \frac{\partial \phi_d}{\partial t}- \text{div}\left(\frac{\phi_d(1-\phi_v-\phi_d)^2}{2M_d(1+\phi_a)^2}\mathbf{T}\nabla \mu \right)=\frac{\Gamma_d}{\gamma} ,\\
\displaystyle \mu=\psi'(\phi_T)-\epsilon^2 \Delta \phi_T,\\
\displaystyle \frac{\partial \phi_a}{\partial t}- \frac{1}{2M_a}\text{div}\left(\mathbf{T}\nabla \phi_a -\chi_a \phi_a \mathbf{T}\nabla c\right)=\frac{\Gamma_a}{\gamma} ,\\
\displaystyle \frac{\partial n}{\partial t}- \alpha_n\text{div}\left(\mathbf{D}_n\nabla n\right)-\alpha_n\chi_v \phi_v= S_n,\\\displaystyle 
\frac{\partial c}{\partial t}- \alpha_c\text{div}\left(\mathbf{D}_c\nabla c\right)-\alpha_c\chi_a \phi_a= S_c,
\end{cases}
\end{equation}
endowed with the homogeneous Neumann boundary conditions
\begin{align}
\label{eqn:19bis}
& b_v\mathbf{T}\nabla \mu \cdot \boldsymbol{\nu}|_{\partial \Omega}=b_d\mathbf{T}\nabla \mu \cdot \boldsymbol{\nu}|_{\partial \Omega}=\left(\mathbf{T}\nabla \phi_a -\chi_a \phi_a \mathbf{T}\nabla c\right)\cdot \boldsymbol{\nu}|_{\partial \Omega}=\nabla \phi_T \cdot \boldsymbol{\nu}|_{\partial \Omega}\\
& \notag=\mathbf{D}_n\nabla n \cdot \boldsymbol{\nu}|_{\partial \Omega}=\mathbf{D}_c\nabla c \cdot \boldsymbol{\nu}|_{\partial \Omega}=0.
\end{align}
System \eqref{eqn:19} is a coupled system of degenerate Cahn--Hilliard with chemotaxis and Keller--Segel equations for the tumor cells and the endothelial cells dynamics respectively. The tensors $\alpha_n\mathbf{D}_n$ and $\alpha_c\mathbf{D}_c$ are related to the diffusion tensor $\mathbf{D}$ of water in the brain tissues, while the tensor of preferential direction $\mathbf{T}$ is related to the brain tissues tractography.

In order to complete the description of the system dynamics, we need to specify biologically meaningful forms for the source terms $\Gamma_v, \Gamma_d, \Gamma_a, S_n, S_c$ in \eqref{eqn:19}. We choose the following forms for the tumor cells growth terms:

\begin{align}
&\frac{\Gamma_v}{\gamma}=\nu \phi_v \max(0,n-\delta_n)(1-\phi_v-\phi_d-\phi_a)-\nu_d\phi_v \max(0,\delta_n-n)-(k_1+k_{T,1})\phi_v,\\
& \notag \frac{\Gamma_d}{\gamma}=k_1\phi_v+\nu_d\phi_v \max(0,\delta_n-n)-(k_2+K_{T,2})\phi_d.
\end{align}
The viable tumor cells proliferate, with a rate $\nu$, as long as the nutrient concentration is above the hypoxia threshold $\delta_n$ and the cells saturation is not reached (contact inhibition), according to a logistic source term. When the nutrient concentration is below the hypoxia threshold, i.e. in hypoxia conditions, the viable cells die by necrosis at a rate $\nu_d$ and proportionally to their concentration. The dead viable cells are transferred to the necrotic phase. The terms $k_1\phi_v$ and $k_2 \phi_v$ indicate clearence terms coming from apoptosis, while $k_{T,1}\phi_v$ and $k_{T,2}\phi_d$ represent possible effects of therapy.
For what concerns the source terms of the chemical species, we choose 
\begin{align}
&S_n=\left(V_n\left(1-H_r(\phi_v+\phi_d)\right)+V_TH_r(\phi_v+\phi_d)(1-\phi_v-\phi_d)\right)(\bar{n}-n)+V_{an}\phi_a(\bar{n}-n)-\delta_v\phi_vn,\\
& \notag S_c=V_c\phi_v\max(0,\delta_n-n)(\bar{c}-c)-\delta_a\phi_ac.
\end{align}
$V_n$ is the supply rate of nutrients from the normal vasculature outside of the tumor, while $V_T$ is the supply rate of nutrients from the normal vasculature inside the tumor. Here, $H_r(\cdot)$ is a regularized heaviside function. Their respective values, depending on the density of the vasculature network and on the pressure levels between the capillaries and the surrounding tissues, are generally different (i.e. inhibited value inside the tumor due to compression of blood vessels). The supply rate inside the tumor is also dependent on the concentration of tumor cells via a factor $(1-\phi_v-\phi_d)$, meaning that the proliferation of the tumor cells destroy blood vessels. Moreover, $\bar{n}$ is the typical nutrient concentration inside the capillaries. The nutrients are released from the capillaries as long as $n<\bar{n}$. The nutrient supply from the tumor induced vasculature is expressed as the product of a specific supply rate $V_{an}$ with the concentration of new vasculature $\phi_a$ and the factor $\bar{n}-n$. Finally, nutrient are consumed by the viable cells at a rate $\delta_v$. For what concerns the source term of the angiogenetic factor, it is released by viable cells at a rate $V_c$ when the nutrient concentration is below the hypoxia threshold and the angiogenetic factor concentration is below its saturation level $\bar{c}$. Moreover, it is consumed by endotelial cells at a rate $\delta_a$.

Finally, the source term for the tumor induced vasculature can be expressed as
\begin{equation}
\frac{\Gamma_a}{\gamma}=\left(1-H_r(\phi_v+\phi_d)\right)\left(V_a\left(1-\phi_v-\phi_d-\phi_a\right)\max(0,c-\delta_c)-k_3\phi_a\right).
\end{equation}
New vessels form at a rate $V_a$ outside of the tumor (see the introduction for details) where the concentration of the angiogenetic factor is greater than a proliferation threshold $\delta_c$ and proportionally to the density of the normal vasculature (which is assumed to be proportional to $\phi_l$), and dissolve by budding and anastomosis at a rate $k_3$.

Reintroducing in \eqref{eqn:19} the parameters $\Pi$ and $k_a$, changing the variables $\tilde{n}=\frac{n}{\bar{n}}$ and $\tilde{c}=\frac{c}{\bar{c}}$ and dividing the fifth and sixth equations of \eqref{eqn:19} by $\bar{n}$ and $\bar{c}$ respectively, we get (omitting the tilde sign)
\begin{equation}
\label{eqn:19tris}
\begin{cases}
\displaystyle
\frac{\partial \phi_v}{\partial t}- \text{div}\left(\frac{\phi_v(1-\phi_v-\phi_d)^2}{L_v(1+\phi_a)^2}\mathbf{T}\nabla \mu \right)+\text{div}\left(h_v\frac{\phi_v(1-\phi_v-\phi_d)^2}{(1+\phi_a)^2}\mathbf{T}\nabla n\right)=\\
 \nu \phi_v \max(0,n-\delta_n)(1-\phi_v-\phi_d-\phi_a)-\nu_d\phi_v \max(0,\delta_n-n)-k_1\phi_v,\\ \\
\displaystyle \frac{\partial \phi_d}{\partial t}- \text{div}\left(\frac{\phi_d(1-\phi_v-\phi_d)^2}{L_d(1+\phi_a)^2}\mathbf{T}\nabla \mu \right)=k_1\phi_v+\nu_d\phi_v \max(0,\delta_n-n)-k_2\phi_d ,\\ \\
\displaystyle \mu=\Pi \psi'(\phi_T)- \Pi \epsilon^2 \Delta \phi_T,\\ \\
\displaystyle \frac{\partial \phi_a}{\partial t}- \frac{1}{L_a}\text{div}\left(\mathbf{T}\nabla \phi_a\right) + h_a \text{div}\left(\phi_a \mathbf{T}\nabla c\right)=\\
\left(1-H_r(\phi_v+\phi_d)\right)\left(V_a\left(1-\phi_v-\phi_d-\phi_a\right)\max(0,c-\delta_c)-k_3\phi_a\right) ,\\ \\
\displaystyle \frac{\partial n}{\partial t}- d_n\text{div}\left(\mathbf{D}_n\nabla n\right)-l_{nv} \phi_v= \\
\left(V_n\left(1-H_r(\phi_v+\phi_d)\right)+V_TH_r(\phi_v+\phi_d)(1-\phi_v-\phi_d)\right)(1-n)+V_{an}\phi_a(1-n)-\delta_v\phi_vn,\\ \\ 
\displaystyle \frac{\partial c}{\partial t}- d_c\text{div}\left(\mathbf{D}_c\nabla c\right)-l_{ca} \phi_a= V_c\phi_v\max(0,\delta_n-n)(1-c)-\delta_a\phi_ac,
\end{cases}
\end{equation}
where $L_v=2M_v$, $L_d=2M_d$, $L_a=\frac{2M_a}{\Pi k_a}$, $h_v=\frac{\Pi \chi_v \bar{n}}{2M_v}$, $h_a=\frac{\Pi \chi_a \bar{c}}{2M_a}$, $d_n=\alpha_n \Pi$, $l_{nv}=\frac{\alpha_n\Pi \chi_v}{\bar{n}}$, $d_c=\alpha_c \Pi$, $l_{ca}=\frac{\alpha_c\Pi \chi_a}{\bar{c}}$. We also redefined the coefficients $\nu=\nu \bar{n}$, $\nu_d=\nu_d \bar{n}$, $\delta_n=\frac{\delta_n}{\bar{n}}$, $V_a=V_a \bar{c}$, $k_3=k_3 \bar{c}$, $\delta_c=\frac{\delta_c}{\bar{c}}$, $V_c=V_c\bar{n}$. 

\section{Matherial and methods}
In this section we describe the data preprocessing workflow to extract reliable quantitative data from the MR images, together with the procedures to generate the patient-specific mesh and the additional meshes containing the values of the components of the water diffusion tensor $\mathbf{D}$, of the tensor of preferential directions $\mathbf{T}$ and of the relative Cerebral Blood Volume (rCBV), starting from the preprocessed MRI, DTI and Dynamic Susceptibility Contrast (DSC) perfusion medical images. 
We then introduce the finite element and time approximation of the continuous problem introduced in Section \ref{sec:theory}, which is implemented to perform numerical simulations of the model on patient--specific test cases.

\subsection{Data preprocessing and mesh generation}
\label{sec:preprocessing}
We integrate clinical data from MRI, DTI and DSC--MRI neuro images of a patient affected by GBM,
 provided by \textit{Fondazione IRCCS Policlinico San Matteo}, into a numerical virtual environment providing a computational mesh of the brain domain and additional meshes containing the values of the independent components of the tensors $\mathbf{D}$ and $\mathbf{T}$ and of the rCBV map.

A mesh generation pipeline is constructed using the analysis tools provided by the FSL \cite{fsl} and FreeSurfer software libraries for image processing and brain tissues segmentation and by the VMTK software library \cite{vmtk} for mesh generation. In particular, starting from a T1--weighted MR image at $1$ mm x $1$ mm x $1$ mm spatial resolution, which provides the structural anatomy of the patient’s brain, the following preprocessing and computational steps are performed:
\begin{itemize}
\item brain extraction by intensity thresholding \cite{bet}, in order to remove non--brain tissues, and bias--field correction \cite{nucorrect};
\item segmentation of the brain tissues and the
background using the FAST algorithm \cite{fast}, based on a hidden Markov random field model and an associated Expectation-Maximization algorithm to estimate the segmentation maps;
\item extraction of a polygonal mesh of the isosurface representing  the external brain boundary
from the segmentation maps, using the marching cubes algorithm \cite{marchingcubes}. Further application of surface smoothing, using Taubin’s algorithm \cite{taubin}, and refinement steps to the isosurface mesh are implemented.
\item Generation of  the 3D mesh using a constrained Delaunay
tetrahedralization of the brain domain defined by its boundary using the
TetGen library \cite{tetgen}, with proper refinement by smooth sizing functions in the
area surrounding the tumour centre.
\end{itemize}
For what concerns the segmentation of the different tumor tissues, we use the nn--Unet deep learning segmentation tool \cite{nnunet}, which is a recently developed state of the art multi task segmentation tool, based on the classical encoding--decoding Unet architecture. When solving the brain tumor segmentation task, starting from a set of T1--weighted, postcontrast T1-weighted, T2--weighted and FLAIR sequences (all sequences brain extracted and co-registered), it generates segmentation labels for the enhancing part of the tumor core, the necrotic part of the tumor core and the peritumoral edema.

A diffusion reconstruction pipeline is further constructed using the MRtrix3 \cite{mrtrix} software library. In particular, starting from raw diffusion data from a DTI sequence comprising a set of single shell diffusion--weighted images at $2.16$ mm x $2.16$ mm x $3.9$ mm spatial resolution with anterior--posterior phase encoding direction and $64$ diffusion--sensitizing directions, the following preprocessing and reconstruction steps are performed:

\begin{itemize}
\item data denoising and noise map estimation by exploiting data redundancy in the PCA domain using random matrix theory \cite{denoise};
\item Gibbs ringing artefacts removal using the method of local subvoxel-shifts \cite{gibbs};
\item eddy current distortion correction by non--linear registration, computed in the phase encoding direction, of the diffusion images onto a single diffusion image \cite{eddy};
\item brain masking, obtained by block-matching rigid registration of the brain mask extracted from the T1--weighted image;
\item non Gaussian diffusion kurtosis estimation of the tensor $\mathbf{D}$ of water diffusion in each voxel using iteratively reweighted linear least squares estimator \cite{tensor2}. This tensor estimation methods, which does not relate on the hypothesis that the diffusion of water molecules follows a Gaussian distribution, is more robust than standard tensor estimation methods based on the Gaussian assumption, since the latter is incorrect for complex biological tissues with cell membranes that create compartments and barriers to diffusion.
\item Projection of the estimated tensor onto the
space of the T1-weighted image by affine block-matching registration.
\end{itemize}
We generate six additional meshes, each one with a piece--wise constant field associated to one independent component of the
tensor $\mathbf{D}$, by assigning to each cell the value of the tensor component of the voxel containing the cell barycentre. 
At the same time, additional six meshes are also associated with each
independent component of the tensor of preferential directions $\mathbf{T}$. This tensor is derived starting from a calculation of the fiber orientation density function from the DTI data using spherical deconvolution \cite{deconv} and further extracting the peaks directions of the spherical harmonic functions in each voxel . The latter method is more robust to derive the fiber orientations inside a voxel  than the methods (both single and  multi compartimental) which use the eigenvectors of the diffusion tensor $\mathbf{D}$. Indeed, this latter method typically fails in regions containing several fiber populations with distinct orientations \cite{deconv}. 

Finally, starting from a 2--d DSC--MRI sequence  at $1.3$ mm x $1.3$ mm x $5$ mm spatial resolution, with echo time $47$ ms , repetition time $1.83$ s, flip angle $90^{\circ}$ and sampling period $90$ s, we implement the preprocessing steps in \cite{perfusion1} and the correction procedure for contrast agent leakage described in \cite{perfusion2} to obtain robust estimate of the rCBV map even in presence of disrupted blood--brain barrier in the peritumoral area. In particular, the following preprocessing and computational steps are performed:
\begin{itemize}
\item spatial smoothing via Gaussian filtering and bias--field correction \cite{nucorrect} of each time frame independently;
\item computation of the tissue and of the arterial concentration--time courses, with a linear and a quadratic dependence on the susceptibility--induced signal drop respectively, after the identification of the baseline signal and the Arterial Input Function in the contralateral Middle Cerebral Artery, as indicated in \cite{perfusion1};
\item quantification of the rCBV map from the time integrals of the concentration--time courses during the passage of the first bolus of contrast agent;
\item corrections to the rCBV values for contrast agent leakage through the estimation of the susceptibility--induced signal drop, via a least--square fit, in terms of a whole--brain average component in nonenhancing voxels and its time integral, as derived in \cite{perfusion2}. 
\end{itemize}

\subsection{Finite element approximation}
We generalize and adapt the Euler semi--implicit finite element approximation introduced in \cite{agosti2}, based on a dual mixed weak formulation of the degenerate Cahn--Hilliard equation and on a discrete variational inequality to impose the positivity of the discrete solution for the tumor concentration, to solve system \eqref{eqn:19tris}.

Let $\mathcal{T}_{h}$ be a quasi-uniform conforming decomposition of $\Omega$ into tetrahedra $K$, 
and let us introduce the following finite element spaces:
\begin{align}
\notag & S_{h} := \{\chi \in C(\bar{\Omega}):\chi |_{K}\in \mathbb{P}^{1}(K) \; \forall K\in \mathcal{T}_{h}\}\subset H^{1}(\Omega),\\
\notag & S_{h}^+ := \{\chi \in S^{h}: \chi \geq 0\; \rm in \, \Omega\}
\end{align}
where $\mathbb{P}_{1}(K)$ indicates the space of polynomials of total order one on $K$.

Let $J$ be the set of nodes of $\mathcal{T}_{h}$ and $\{\mathbf{x}_j\}_{j\in J}$ be the set of their coordinates. Moreover, let $\{\varphi_j\}_{j\in J}$ be the Lagrangian basis functions associated with each node $j\in J$.
Denoting by $\pi^h:C(\bar{\Omega})\rightarrow S^h$ the standard Lagrangian interpolation operator we define the lumped scalar product as
\begin{equation}
\label{eqn:lump}
(\eta_1,\eta_2)^h=\int_{\Omega}\pi^h(\eta_1(\mathbf{x})\eta_2(\mathbf{x}))d\mathbf{x}\equiv \sum_{j\in J}(1,\chi_j)\eta_1(\mathbf{x}_j)\eta_2(\mathbf{x}_j),
\end{equation}
for all $\eta_1,\eta_2\in C(\bar{\Omega})$.

We set $\Delta t = T/N$ for a $N \in \mathbb{N}$ and $t_{n}=n\Delta t$, {$n=0,...,N$}, where $[0,T]$ is the time span of the system dynamics.
Starting from data $\phi_{v}^0, \phi_{d}^0, \phi_{a}^0, n^{0}, c^0,  \in \left[H^{2}(\Omega)\right]^5$ and $\phi_{h,v}^{0}=\pi^{h}\phi_{v}^0$, $\phi_{h,d}^{0}=\pi^{h}\phi_{d}^0$, $\phi_{h,a}^{0}=\pi^{h}\phi_{a}^0$, $n_{h}^{0}=\pi^{h}n^0$, $c_{h}^{0}=\pi^{h}c^0$, with $0\leq \phi_{h,v}^{0}, \phi_{h,d}^{0}, \phi_{h,a}^{0}, n_{h}^{0}, c_{h}^{0},<1$, we consider the following fully discretized problem:
for $n=1,\dots,N$, given $(\phi_{h,v}^{n-1}, \phi_{h,d}^{n-1}, \phi_{h,a}^{n-1}, n_h^{n-1}, c_h^{n-1}) \in \left[S_h^+\right]^5$, find $(\phi_{h,v}^{n},\Sigma_{h,v}^{n}, \phi_{h,d}^{n},\Sigma_{h,d}^{n}, \phi_{h,a}^{n}, n_h^{n}, c_h^{n}) \in S_h^+ \times S_h \times S_h^+ \times S_h \times S_h^+ \times S_h^+ \times S_h^+$ such that, $\forall (v_h,w_h,q_h, r_h, s_h, t_h, u_h) \in S_h \times S_h^+ \times S_h \times S_h^+ \times S_h \times S_h \times S_h$,
\begin{equation}
\label{eqn:discr}
\begin{cases}
\displaystyle 
\biggl ( \frac{\phi_{h,v}^{n} -\phi_{h,v}^{n-1}}{\Delta t} , v_h \biggr )^h + \frac{1}{L_v} \left( \frac{\phi_{h,v}^{n-1}(1-\phi_{h,v}^{n-1}-\phi_{h,d}^{n-1})^2}{(1+\phi_{h,a}^{n-1})^2} \mathbf{T} \nabla \Sigma_{h,v}^{n} , \nabla v_h \right)   \\
= \nu (  \phi_{h,v}^{n-1} \max(0,n_h^{n}-\delta_n)(1-\phi_{h,v}^{n-1}-\phi_{h,d}^{n-1}-\phi_{h,a}^{n-1}) , v_h )^h   \\
 - \nu_d (\phi_{h,v}^{n-1}\max(0,\delta_n-n_h^n), v_h)^h +h_v \left( \frac{\phi_{h,v}^{n-1}(1-\phi_{h,v}^{n-1}-\phi_{h,d}^{n-1})^2}{(1+\phi_{h,a}^{n-1})^2}\mathbf{T} \nabla n_h^{n} , \nabla v_h  \right)\\
 -((k_1+k_{T,1})\phi_{h,v}^{n-1},v_h)^h,
\\ \\
\displaystyle 
\biggl ( \frac{\phi_{h,d}^{n} -\phi_{h,d}^{n-1}}{\Delta t} , q_h \biggr )^h + \frac{1}{L_d} \left( \frac{\phi_{h,d}^{n-1}(1-\phi_{h,v}^{n-1}-\phi_{h,d}^{n-1})^2}{(1+\phi_{h,a}^{n-1})^2} \mathbf{T} \nabla \Sigma_{h,d}^{n} , \nabla q_h \right)   \\
 = \nu_d (\phi_{h,v}^{n-1}\max(0,\delta-n_h^n), q_h)^h +(k_1\phi_{h,v}^{n-1},q_h)^h-((k_2+k_{T,2})\phi_{h,d}^{n-1},q_h)^h ,
\\ \\
\displaystyle  \Pi \epsilon^2\left(\nabla \left(\phi_{h,v}^n+\phi_{h,d}^n\right),\nabla(w_h-\phi_{h,v}^n)\right)+(\Pi \psi'_1(\phi_{h,v}^n+\phi_{h,d}^n),w_h-\phi_{h,v}^n)^h \geq \\ 
\displaystyle (\Sigma_{h,v}^n-\Pi \psi'_2(\phi_{h,v}^{n-1}+\phi_{h,d}^{n-1}),w_h-\phi_{h,v}^n)^h, \\ \\
\displaystyle \Pi \epsilon^2\left(\nabla \left(\phi_{h,v}^n+\phi_{h,d}^n\right),\nabla(r_h-\phi_{h,d}^n)\right)+(\Pi \psi'_1(\phi_{h,v}^n+\phi_{h,d}^n),r_h-\phi_{h,d}^n)^h \geq \\
\displaystyle (\Sigma_{h,d}^n-\Pi \psi'_2(\phi_{h,v}^{n-1}+\phi_{h,d}^{n-1}),r_h-\phi_{h,d}^n)^h, \\ \\
\biggl ( \frac{\phi_{h,a}^{n} -\phi_{h,a}^{n-1}}{\Delta t} , s_h \biggr )^h + \frac{1}{L_a} \left( \mathbf{T} \nabla \phi_{h,a}^{n} , \nabla s_h \right) -h_a \left( \phi_{h,a}^{n-1}\mathbf{T} \nabla c_h^{n} , \nabla v_h  \right) \\
= \left(\left(1-H_r(\phi_{h,v}^{n-1}+\phi_{h,d}^{n-1})\right)\left(V_a\left(1-\phi_{h,v}^{n-1}-\phi_{h,d}^{n-1}-\phi_{h,a}^{n-1}\right)\max(0,c_h^n-\delta_c)-k_3\phi_{h,a}^{n-1}\right),s_h\right)^h ,\\ \\
\displaystyle  \biggl ( \frac{n_h^{n} -n_h^{n-1}}{\Delta t} , t_h \biggr )^h + b_n( \mathbf{D} \nabla n_h^{n} , \nabla t_h ) = \left(l_{nv} \phi_{h,v}^{n-1},t_h\right)^h\\
+ \biggl(\left(V_n\left(1-H_r(\phi_{h,v}^{n-1}+\phi_{h,d}^{n-1})\right)+V_TH_r(\phi_{h,v}^{n-1}+\phi_{h,d}^{n-1})(1-\phi_{h,v}^{n-1}-\phi_{h,d}^{n-1})\right)(1-n_h^n)\\
+V_{an}\phi_{a,h}^{n-1}(1-n_h^n)-\delta_v\phi_{h,v}^{n-1}n_h^n,t_h\biggr)^h,\\ \\
\displaystyle  \biggl ( \frac{c_h^{n} -c_h^{n-1}}{\Delta t} , u_h \biggr )^h + b_c( \mathbf{D} \nabla c_h^{n} , \nabla u_h ) = \left(l_{ca} \phi_{h,a}^{n-1},u_h\right)^h+\\\left(\left(V_c\phi_{h,v}^{n-1}\max(0,\delta_n-n_h^n)(1-c_h^n)-\delta_a\phi_{h,a}^{n-1}c_h^n\right),u_h\right)^h,
\end{cases}
\end{equation}
where $(\cdot,\cdot)$ denotes the standard $L^2$ inner product over $\Omega$, and the following convex splitting of the potential is considered:
\vspace{-0.3cm}
\begin{equation}
\label{eqn:convsplit}
\psi_1(\phi)=-(1-\bar{\phi})\ln(1-\phi), \quad \psi_2(\phi)=-\biggl[\frac{1}{3}\phi^3+\frac{1-\bar{\phi}}{2}\phi^2+(1-\bar{\phi})\phi \biggr].
\end{equation}
In the fifth and sixth equations of \eqref{eqn:discr} we have supposed that the tensors $d_n \mathbf{D}_n$ and $d_c \mathbf{D}_c$ are proportional to the diffusion tensor of water molecules $\mathbf{D}$, with proportionality constants $b_n$ and $b_c$ respectively. 

The finite element approximation \eqref{eqn:discr} has the form of a coupled system of discrete variational inequalities, where the positivity of the discrete
solutions $\phi_{h,v}$ and $\phi_{h,n}$ is enforced as a constraint, being independently projected onto the space with positive values $S_h^+$ . 
\begin{rem}
\label{rem:1}
The positivity of the discrete solutions $\phi_{h,a}, n_h, c_h$ of \eqref{eqn:discr} is guaranteed only for sufficiently small values of $\Delta t$. The design of a proper well posed, convergent and gradient--stable approximation of system \eqref{eqn:19}, where the Keller--Segel equations are discretized  in a consistent way with the gradient flow structure in \eqref{eqn:11} and where the positivity of $\phi_{h,a}, n_h, c_h$ is enforced by specific approximations schemes for the corresponding second order parabolic equations, will be studied as a future development. For the purpose of the present work, we check at each time step of the numerical simulations if the discrete solutions $\phi_{h,a}, n_h, c_h$ satisfy the positivity constraints and the saturation constraint, and in case of violation of these constraints we reduce the time step in an adaptive way.
\end{rem}
We observe that the equation for $n_h^n$ in the semi--implicit scheme in \eqref{eqn:discr} is decoupled from the other equations, and can be solved independently. Once we obtain a solution $n_h^n$, we can solve for $c_h^n$ and for $\phi_{h,a}^n$ sequentially. Given $n_h^n, c_h^n, \phi_{h,a}^n$, a solution for the coupled system of equations and inequalities for $\phi_{h,v}^n, \Sigma_{h,v}^n, \phi_{h,d}^n, \Sigma_{h,d}^n$ can be obtained.
\newline
The lumped mass approximation of the $L^2$ scalar product is introduced in  \eqref{eqn:discr} in order for the discrete solutions $\phi_{h,v}^n$ and $\phi_{h,d}^n$
to be able to select compactly supported solutions with a moving free boundary from the unphysical ones with fixed support (see \cite{agosti2,agosti3} for discussions).
Note also that the convex part of the cellular potential, depending on both $\phi_{h,v}^{n}$ and $\phi_{h,d}^{n}$, is treated implicitly in time, whereas the concave part is treated explicitly. 
  
In order to solve the discrete system \eqref{eqn:discr} for every time step $n$, we generalize the iterative procedure introduced in \cite{agosti2}. Let us introduce two families of subdivisions of the set $J$ in the following way: we define a set of \textit{viable passive nodes} $J_{0,v}(\phi_{h,v}^{n-1})\subset J$ such that $ (\phi_{h,v}^{n-1},\varphi_j)=0$ for $j\in J_{0,v}(\phi_{h,v}^{n-1})$ and a set of \textit{viable active nodes} $J_{+,v}(\phi_{h,v}^{n-1})=J \setminus J_{0,v}(\phi_{h,v}^{n-1})$; we also define a set of \textit{necrotic passive nodes} $J_{0,d}(\phi_{h,d}^{n-1})\subset J$ such that $ (\phi_{h,d}^{n-1},\varphi_j)=0$ for $j\in J_{0,d}(\phi_{h,d}^{n-1})$ and a set of \textit{necrotic active nodes} $J_{+,d}(\phi_{h,d}^{n-1})=J \setminus J_{0,d}(\phi_{h,d}^{n-1})$. A viable (resp. necrotic) passive node is thus characterized by the fact that $\phi_{h,v}^{n-1}\equiv 0$ (resp. $\phi_{h,d}^{n-1}\equiv 0$) on the support of the basis function associated to it.
We formulate the following algorithm: 
\begin{algorithmic}
\Require $\mu>0$ (a relaxation parameter), $\phi_{h,v}^{n-1},\Sigma_{h,v}^{n-1},\phi_{h,d}^{n-1},\Sigma_{h,d}^{n-1},\phi_{h,a}^{n-1},n_{h}^{n-1},c_{h}^{n-1}$;\\
\textbf{Step 0} Find $(\phi_{h,a}^{n},n_h^n,c_h^n)\in S_h^+ \times S_h^+ \times S_h^+$, $\forall (s_h,t_h,u_h) \in S_h \times S_h \times S_h$, such that:
\begin{equation*}
\begin{cases}
\biggl ( \frac{\phi_{h,a}^{n} -\phi_{h,a}^{n-1}}{\Delta t} , s_h \biggr )^h + \frac{1}{L_a} \left( \mathbf{T} \nabla \phi_{h,a}^{n} , \nabla s_h \right) +h_a  \left( \phi_{h,a}^{n-1}\mathbf{T} \nabla c_h^{n} , \nabla v_h  \right) \\
= \left(\left(1-H_r(\phi_{h,v}^{n-1}+\phi_{h,d}^{n-1})\right)\left(V_a\left(1-\phi_{h,v}^{n-1}-\phi_{h,d}^{n-1}-\phi_{h,a}^{n-1}\right)\max(0,c_h^n-\delta_c)-k_3\phi_{h,a}^{n-1}\right),s_h\right)^h ,\\ \\
\displaystyle  \biggl ( \frac{n_h^{n} -n_h^{n-1}}{\Delta t} , t_h \biggr )^h + b_n( \mathbf{D} \nabla n_h^{n} , \nabla t_h ) = (l_{nv} \phi_{h,v}^{n-1},t_h)^h\\
+ \biggl(\left(V_n\left(1-H_r(\phi_{h,v}^{n-1}+\phi_{h,d}^{n-1})\right)+V_TH_r(\phi_{h,v}^{n-1}+\phi_{h,d}^{n-1})(1-\phi_{h,v}^{n-1}-\phi_{h,d}^{n-1})\right)(1-n_h^n)\\
+V_{an}\phi_{a,h}^{n-1}(1-n_h^n)-\delta_v\phi_{h,v}^{n-1}n_h^n,t_h\biggr)^h,\\ \\
\displaystyle  \biggl ( \frac{c_h^{n} -c_h^{n-1}}{\Delta t} , u_h \biggr )^h + b_c( \mathbf{D} \nabla c_h^{n} , \nabla u_h ) = (l_{ca} \phi_{h,a}^{n-1},u_h)^h+\\\left(\left(V_c\phi_{h,v}^{n-1}\max(0,\delta_n-n_h^n)(1-c_h^n)-\delta_a\phi_{h,a}^{n-1}c_h^n\right),u_h\right)^h,
\end{cases}
\end{equation*} \\
\For{$k\geq 0$}  \\
\textbf{Initialization}\[\phi_{h,v}^{n,0}=\phi_{h,v}^{n-1},\Sigma_{h,v}^{n,0}=\Sigma_{h,v}^{n-1},\phi_{h,vd}^{n,0}=\phi_{h,d}^{n-1},\Sigma_{h,d}^{n,0}=\Sigma_{h,d}^{n-1};\]
\textbf{Step 1} Find $(z_v^{n,k},z_d^{n,k})\in S_h\times S_h$ such that $\forall (v_h,w_h) \in S_h\times S_h$:
\[(z_v^{n,k},v_h)^h=(\phi_{h,v}^{n,k},v_h)^h-\mu [\Pi \epsilon^2 \left(\nabla \left(\phi_{h,v}^{n,k}+\phi_{h,d}^{n,k}\right),\nabla v_h\right)+(\Pi \psi_{2}'(\phi_{h,v}^{n-1}+\phi_{h,d}^{n-1})-\Sigma_{h,v}^{n,k},w_h)^h],\]
\[(z_d^{n,k},w_h)^h=(\phi_{h,d}^{n,k},w_h)^h-\mu [\Pi\epsilon^2 \left(\nabla \left(\phi_{h,v}^{n,k}+\phi_{h,d}^{n,k}\right),\nabla w_h\right)+(\Pi \psi_{2}'(\phi_{h,v}^{n-1}+\phi_{h,d}^{n-1})-\Sigma_{h,d}^{n,k},w_h)^h];\]
\textbf{Step 2} Find $(\phi_{h,v}^{n,k+1/2},\phi_{h,d}^{n,k+1/2})\in S_h^+\times S_h^+$, $\forall r,s\geq 0$, such that:

\For{$j\in J$}
\If {$j \in J_{0,v}(\phi_{h,v}^{n-1})$}
    \State $\phi_{h,v}^{n,k+1/2}(\mathbf{x}_j)\gets \phi_{h,v}^{n-1}(\mathbf{x}_j)$
\Else
        \begin{equation}
        \label{eqn:varinequalv}
        (\phi_{h,v}^{n,k+1/2}(\mathbf{x}_j)+\mu \Pi \psi_{1}'(\phi_{h,v}^{n,k+1/2}+\phi_{h,d}^{n,k})(\mathbf{x}_j)-z_v^{n,k}(\mathbf{x}_j),r-\phi_{h,v}^{n,k+1/2}(\mathbf{x}_j))\geq 0              \end{equation}
        \EndIf  \EndFor 
\For{$j\in J$}
\If {$j \in J_{0,d}(\phi_{h,d}^{n-1})$}
    \State $\phi_{h,d}^{n,k+1/2}(\mathbf{x}_j)\gets \phi_{h,d}^{n-1}(\mathbf{x}_j)+\Delta t\left(\nu_d (\phi_{h,v}^{n-1}\max(0,\delta-n_h^n))+(k_1\phi_{h,v}^{n-1})\right)(\mathbf{x}_j)$
\Else
        \begin{equation}
        \label{eqn:varinequald}
        (\phi_{h,d}^{n,k+1/2}(\mathbf{x}_j)+\mu \Pi \psi_{1}'(\phi_{h,v}^{n,k}+\phi_{h,d}^{n,k+1/2})(\mathbf{x}_j)-z_d^{n,k}(\mathbf{x}_j),s-\phi_{h,s}^{n,k+1/2}(\mathbf{x}_j))\geq 0              \end{equation}
        \EndIf  \EndFor \\
\textbf{Step 3} Find $(\phi_{h,v}^{n,k+1},\Sigma_{h,v}^{n,k+1},\phi_{h,d}^{n,k+1},\Sigma_{h,d}^{n,k+1})\in S_h^+ \times S_h \times S_h^+ \times S_h$, $\forall (v_h,w_h,q_h,r_h) \in S_h \times S_h \times S_h \times S_h$, such that:
\begin{equation*}
\begin{cases}
\displaystyle \frac{1}{\Delta t} ( \phi_{h,v}^{n, k+1},v_h )^h + \frac{1}{L_v}  ( \mathbf{T} \nabla \Sigma_{h,v}^{n,k+1},\nabla v_h )=\frac{1}{\Delta t} ( \phi_{h,v}^{n-1},v_h )^h+h_v \left( \frac{\phi_{h,v}^{n-1}(1-\phi_{h,v}^{n-1}-\phi_{h,d}^{n-1})^2}{(1+\phi_{h,a}^{n-1})^2}\mathbf{T} \nabla n_h^{n} , \nabla v_h  \right) \\ \displaystyle + \nu (  \phi_{h,v}^{n-1} \max(0,n_h^{n}-\delta_n)(1-\phi_{h,v}^{n-1}-\phi_{h,d}^{n-1}-\phi_{h,a}^{n-1}) , v_h )^h  - \nu_d (\phi_{h,v}^{n-1}\max(0,\delta_n-n_h^n), v_h)^h  \\ \displaystyle
 -((k_1+k_{T,1})\phi_{h,v}^{n-1},v_h)^h+ \frac{1}{L_v} \left( \left[1-\frac{\phi_{h,v}^{n-1}(1-\phi_{h,v}^{n-1}-\phi_{h,d}^{n-1})^2}{(1+\phi_{h,a}^{n-1})^2}\right] \mathbf{T} \nabla \Sigma_{h,v}^{n,k} , \nabla v_h \right)\\ \\
\displaystyle ( \phi_{h,v}^{n,k+1},w_h   )^h + \mu \Pi \epsilon^2 \left(\nabla\left(\phi_{h,v}^{n,k+1}+\phi_{h,d}^{n,k+1}\right),\nabla w_h\right) -\mu  (\Sigma_{h,v}^{n,k+1},w_h)^h = \\
\displaystyle ( 2\phi_{h,v}^{n,k+1/2}-z_{h,v}^{n,k}-\mu \Pi \psi'_2(\phi_{h,v}^{n-1}+\phi_{h,v}^{n-1}),w_h  )^h,\\ \\
\displaystyle \frac{1}{\Delta t} ( \phi_{h,d}^{n, k+1},q_h )^h + \frac{1}{L_d}  ( \mathbf{T} \nabla \Sigma_{h,d}^{n,k+1},\nabla q_h )=\frac{1}{\Delta t} ( \phi_{h,d}^{n-1},q_h )^h+ \nu_d (\phi_{h,v}^{n-1}\max(0,\delta-n_h^n), q_h)^h \\
\displaystyle +(k_1\phi_{h,v}^{n-1},q_h)^h - ((k_2+k_{T,2})\phi_{h,d}^{n-1},q_h)^h+ \frac{1}{L_d} \left( \left[1-\frac{\phi_{h,d}^{n-1}(1-\phi_{h,v}^{n-1}-\phi_{h,d}^{n-1})^2}{(1+\phi_{h,a}^{n-1})^2}\right] \mathbf{T} \nabla \Sigma_{h,d}^{n,k} , \nabla q_h \right)\\ \\
\displaystyle ( \phi_{h,d}^{n,k+1},r_h   )^h + \mu \Pi \epsilon^2 \left(\nabla\left(\phi_{h,v}^{n,k+1}+\phi_{h,d}^{n,k+1}\right),\nabla r_h\right) -\mu  (\Sigma_{h,d}^{n,k+1},r_h)^h = \\
\displaystyle ( 2\phi_{h,d}^{n,k+1/2}-z_{h,d}^{n,k}-\mu \Pi \psi'_2(\phi_{h,v}^{n-1}+\phi_{h,v}^{n-1}),r_h  )^h.
\end{cases}
\end{equation*} \\
\If {$||\phi_{h,v}^{n,K+1}-\phi_{h,v}^{n,K}||_{\infty}+||\phi_{h,d}^{n,K+1}-\phi_{h,d}^{n,K}||_{\infty}<10^{-6}$}
    \State $(\phi_{h,v}^n,\Sigma_{h,v}^n,\phi_{h,d}^n,\Sigma_{h,d}^n)\gets (\phi_{h,v}^{n,K+1},\Sigma_{h,v}^{n,K+1},\phi_{h,d}^{n,K+1},\Sigma_{h,d}^{n,K+1})$; \textbf{break}.
\EndIf
\EndFor
\end{algorithmic}
The variational inequalities \eqref{eqn:varinequalv} and \eqref{eqn:varinequald} contain the elliptic terms in the forcing terms $z_v^{n,k}$ and $z_d^{n,k}$ only, which are known terms at the step $k+1/2$, and hence can be solved as projection problems on each active node independently.
In order to solve them, we solve the following projected gradient methods in the index $l$, $l=0, \dots, L$, starting from $\phi_{h,v}^{n,k+1/2,0}=\phi_{h,v}^{n,k}$ and $\phi_{h,d}^{n,k+1/2,0}=\phi_{h,d}^{n,k}$:
\begin{align}
\label{eq:projgrad}
&\phi_{h,v}^{n,k+1/2,l+1}(\mathbf{x}_j)=\max \biggl \{0, \phi_{h,v}^{n,k+1/2,l}(\mathbf{x}_j) -\omega [  \phi_{h,v}^{n,k+1/2,l}(\mathbf{x}_j) +\mu \Pi \psi'_1(\phi_{h,v}^{n,k+1/2,l}+\phi_{h,d}^{n,k})(\mathbf{x}_j)-z_{h,v}^{n,k} ] \biggr \},\\
& \notag \phi_{h,d}^{n,k+1/2,l+1}(\mathbf{x}_j)=\max \biggl \{0, \phi_{h,d}^{n,k+1/2,l}(\mathbf{x}_j) -\omega [  \phi_{h,d}^{n,k+1/2,l}(\mathbf{x}_j) +\mu \Pi \psi'_1(\phi_{h,v}^{n,k}+\phi_{h,d}^{n,k+1/2,l})(\mathbf{x}_j)-z_{h,d}^{n,k} ] \biggr \},
\end{align} 
where $\omega$ is a relaxation parameter. If {$\| \phi_{h,v}^{n,k+1/2,L+1} - \phi_{h,v}^{n,k+1/2,L}  \|_{\infty}+\| \phi_{h,d}^{n,k+1/2,L+1} - \phi_{h,d}^{n,k+1/2,L}  \|_{\infty}<10^{-6}$}, we stop the cycle and set $\phi_{h,v}^{n,k+1/2}=\phi_{h,v}^{n,k+1/2,L+1}, \phi_{h,d}^{n,k+1/2}=\phi_{h,d}^{n,k+1/2,L+1}$. Note that, since the operators acting on $\phi_{h,v}^{n,k+1/2,l}$ and $\phi_{h,d}^{n,k+1/2,l}$ in the square brackets in \eqref{eq:projgrad} are continuous and strictly monotone, the projection maps defined in \eqref{eq:projgrad} have a unique fixed point \cite{temam}.

\section{Results}
In this section we present the results obtained by numerically solving \eqref{eqn:discr} on patient--specific data, generated following the procedures reported in section \ref{sec:preprocessing}.

The numerical algorithm is implemented in the FEniCS computing platform \cite{fenics}, a collection of software libraries to solve partial differential equations with the finite element methods, with high-level Python and C++ interfaces.

As in \cite{agosti2}, we implement an adaptive time step based on the \textit{CFL}--type condition for the viable cells phase, i.e. at each time level n we set $\Delta t ^n< h_{\text{min}}/v_{\text{max},v}^{n-1}$, where $h_{\text{min}}$ is the smallest edge length among the {mesh cells} and $v_{\text{max},v}^{n-1}$ is the maximum on $\Omega$ of the ${\infty}-$norm of the viable cells velocity $\mathbf{v}_v^{n-1}$ at the previous time step, obtained rewriting the first equation of \eqref{eqn:19tris} as $\frac{\partial \phi_v}{\partial t}+\text{div}(\phi_v\mathbf{v}_v)=\frac{\Gamma_v}{\gamma}$, i.e.
\begin{equation}
\label{eqn:vel}
\mathbf{v}_v^{n-1}= h_v \frac{(1-\phi_{h,v}^{n-1}-\phi_{h,d}^{n-1})^2}{(1+\phi_{h,a}^{n-1})^2}\mathbf{T} \nabla n_h^{n-1}-\frac{(1-\phi_{h,v}^{n-1}-\phi_{h,d}^{n-1})^2}{L_v(1+\phi_{h,a}^{n-1})^2}\mathbf{T} \nabla \Sigma_{h,v}^{n-1},
\end{equation}
and
\begin{equation}
\label{eqn:vmax}
\mathbf{v}_{\text{max},v}^{n-1}:=\max_{j\in J}(|\mathbf{v}_{v,x}(\mathbf{x}_j)|+|\mathbf{v}_{v,y}(\mathbf{x}_j)|+|\mathbf{v}_{v,z}(\mathbf{x}_j)|).
\end{equation}
We thus impose
\begin{equation}
\label{eqn:deltat}
\Delta t ^n= \min \biggl (100 \frac{L_v}{\Pi}\epsilon^2, \frac{h_{\text{min}}}{ 2\mathbf{v}_{\text{max},v}^{n-1}} \biggr),
\end{equation}
where $M_v\epsilon^2/\Pi$ is the typical time scale for the spinodal decomposition dynamics of $\phi_v$.
We moreover halve the time step \eqref{eqn:deltat} whenever the concentrations $\phi_{h,a}^n$ $n_h^n$ or $c_h^n$ violates the physical constraints (positivity and mixture saturation). 

We choose the following values of the model parameters in \eqref{eqn:discr} for the tumor cells dynamics, which belong to the range of biological parameters reported in \cite{agosti2,agosti3}: $\Pi=694$ Pa, $L_v=L_d=3900$ (Pa day)/mm$^2$, $h_v=0.14\left(\chi(CSF,GM)+4\chi(WM)\right)$ $mm^2/\text{day}$ (obtained by taking $\bar{n}=0.07$ mMol), $\nu = 0.15$ 1/day, $\nu_d = 0.06$ 1/day,  $\phi_e = 0.389$, $\epsilon = 0.013$ mm, $\delta_n=0.33$, $k_1=k_2=0$, $b_n=1$, $V_n=10^4$ 1/day,  $\delta_v=8640$ 1/day, $l_{nv}=111.42\left(\chi(CSF,GM)+4\chi(WM)\right)$ 1/day (obtained by taking $\alpha_n=0.001$ $mMol^2$/Pa day). Here, $\chi(T)$ is the characteristic function of the tissue T, with T=CerebroSpinal Fluid (CSF), Gray Mater (GM) or White Matter (WM). $b_n$ is taken to be equal to $1$, in the hypothesis that the nutrient (i.e. Oxygen) is perfectly diluted in water. The basic temporal step is thus $\Delta t=100 \cdot L_v\epsilon^2/\Pi = 0.095$ day. 
\newline
For what concerns the endothelial cells dynamics, we choose the following values of the model parameters in \eqref{eqn:discr}, which belong to the range of biological parameters reported in \cite{chaplain,frieboes}: $L_a^{-1}=0.003$ mm$^2$/Pa, $h_a=0.2264$ mm$^2$/Pa, $\delta_c=0.2$, $V_a=4.8$ 1/day, $k_3=0.24$ 1/day, $b_c=0.589$ (obtained by considering a reference diffusion coefficient for the angiogenetic factor of $50.9$ mm$^2$/day, and a reference diffusion coefficient for Oxygen of $86.4$ mm$^2$/day), $l_{ca}=0.73$ 1/day (obtained by taking $\alpha_c=0.001$ $mMol^2$/Pa day and $\bar{c}=0.1$ mMol), $V_c=10^3$ 1/day, $\delta_a=864$ 1/day. We observe that there are no estimates in literature of the parameters $V_c$ and $\delta_a$ by in vitro or in vivo experiments, so we take values with the same order of magnitude to the ones proposed in the numerical simulations in \cite{chaplain}, which are one order of magnitude less than $V_n$ and $\delta_v$.
\newline
We finally need to assign specific values to the parameters $V_T$ and $V_{an}$, which represent the source rates of nutrient from the normal vasculature inside the tumor and from the new vasculature respectively. To the best of our knowledge, there are no empirical estimates for the values of these parameters in literature. In the following, we take $V_{an}=V_T/\bar{\phi}_a$, where $\bar{\phi}_a$ is an estimate for the average blood volume concentration in the brain. This means that, when the concentration of the tumor induced vasculature is the same as the average vasculature concentration in the healthy tissues, i.e. $\phi_a=\bar{\phi}_a$, then the supply rate from the new vasculature is proportional to $V_{an}\phi_a=V_T$. We note that $V_T$ can be related to the value of $V_n$, which has been estimated e.g. in \cite{agosti2,agosti3} in correspondence to a normal vasculature concentration, and we expect the new vasculature to supply the same amount of nutrient as the preexisting vasculature inside the tumor when its density is comparable with the preexisting vasculature density. We then take $\bar{\phi}_a=0.04$, based on the considerations that the average blood volume in an healthy brain is $3.5$ $ml/100 g$ in Grey Matter and $1.7$ $ml/100 g$ in White Matter, while the average density of an healthy brain is $1081$  $g/ml$ \cite{bv1,bv2,bv3}. We finally consider two possible values for $V_T$: $V_T=V_n/2$ (\textbf{Case 1}) and $V_T=V_n/10$ (\textbf{Case 2}). Lowering the value of $V_T$, i.e. the nutrient supply inside the tumor tissue, will trigger and accelerate the angiogenetic process induced by the production of the angiogenetic factor by viable cells in hypoxia condition and its subsequent diffusion in the peritumoral area.
A sensitivity analysis and optimization of the model parameters, based on longitudinal patient--specific data, will be studied in a forthcoming paper.

In the following, starting from the neuroimaging data from a patient affected by GlioBlastoma Multiforme acquired at different temporal stages, i.e. before surgery, after surgery and after therapy, we show the numerical results about three test cases. In particular, in \textbf{Test case 1}, based on data before surgery, we consider the evolution of a small early stage tumor, without a necrotic core, artificially placed in a healthy region of the brain. In \textbf{Test case 2}, based on data before surgery, we consider the evolution of the real tumor, as delineated through the automatic segmentation process described in section \ref{sec:preprocessing}, which comprises both a viable and a necrotic tumor region. The simulated time span for both Test case 1 and Test Case 2 is $60$ days. In both these studies, we consider the parameters set with high nutrient supply inside the tumor tissues (\textbf{Case 1}) and with low nutrient supply inside the tumor tissues (\textbf{Case 2}), in order to observe the sensitivity of the angiogenetic process to the availability of nutrients inside the tumor. Finally, in \textbf{Test case 3}, based on data after surgery, we consider the evolution of residual tumor particles randomly placed on the boundary of the resection area, under the application of therapy according to the standard Stupp protocol \cite{stupp}, and we compare the simulation results with data acquired during a control after therapy, seven months after the surgery event. For the considered subject, no perfusion data were acquired before and after surgery, while the DSC--MRI data are available only in the after therapy event. Therefore we compare the simulation results for both the tumor and the tumor induced vasculature distributions with the available structural and perfusion data in \textbf{Test case 3}.
\begin{figure}[ht!]
\includegraphics[width=0.9\linewidth]
{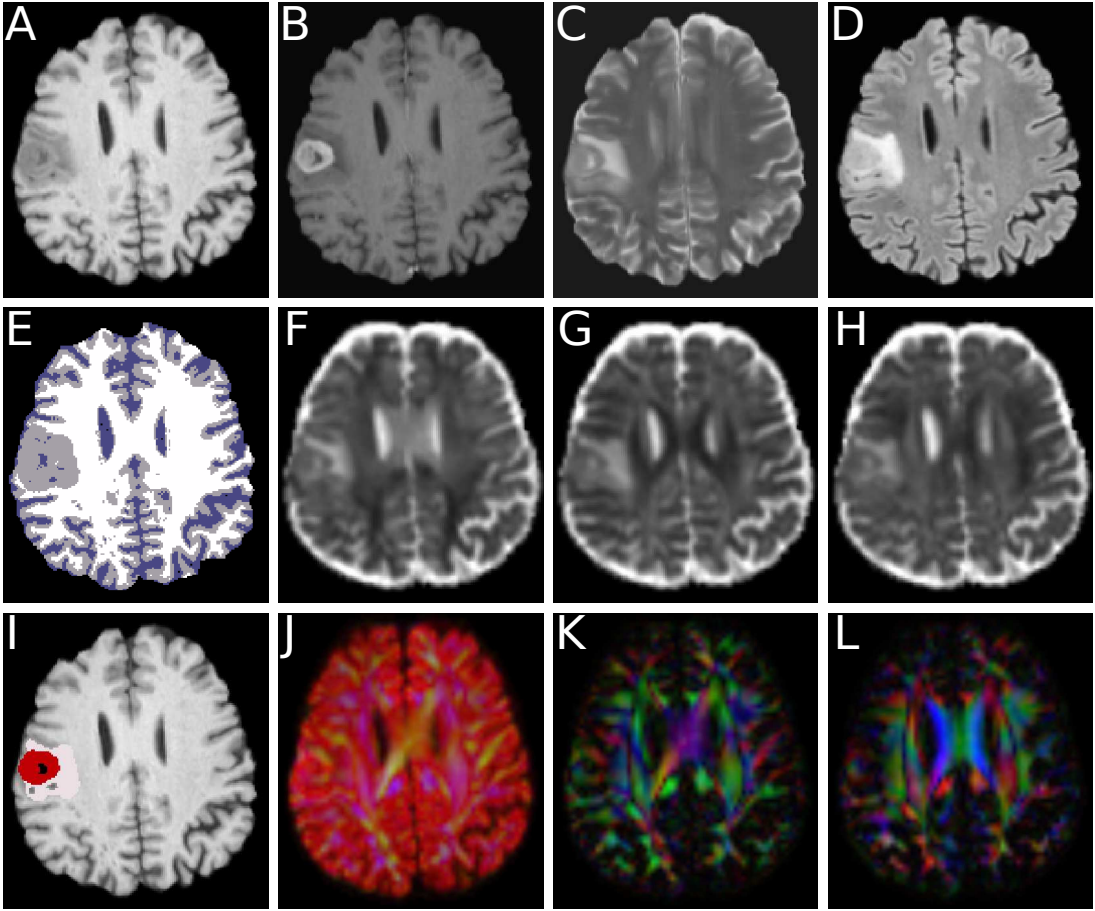}
\centering
\caption{Axial slice from the T1--weighted MRI (A), the postcontrast T1--weighted MRI (B), the T2--weighted MRI (C), the FLAIR MRI (D), the segmentation map of the brain tissues into WM (white color), GM (gray color) and CSF (blue color) (E), the $xx$ component of the diffusion tensor $\mathbf{D}$ (F), the $yy$ component of the diffusion tensor $\mathbf{D}$ (G), the $zz$ component of the diffusion tensor $\mathbf{D}$ (H),  the tumor tissues segmentation into a viable component (brown color), a necrotic component (black color) and edema (pink color) (I), the coefficient of the first spherical harmonics decomposition of the diffusion tensor, in RGB color scale (J), the coefficient of the second spherical harmonics decomposition of the diffusion tensor, in RGB color scale (K), the coefficient of the third spherical harmonics decomposition of the diffusion tensor, in RGB color scale (L).}
\label{fig:1}
\end{figure}
In Figure \ref{fig:1} we show a representative illustration of the (preprocessed) neuroimaging data from the considered subject acquired before surgery, reported on a selected axial slice. In the first row of the figure, the T1--weighted MRI (A), the postcontrast T1--weighted MRI (B), the T2--weighted MRI (C) and the FLAIR MRI (D) are shown. In the second row, the segmentation map of the brain tissues into WM (white color), GM (gray color) and CSF (blue color) (E), the $xx$ component of the diffusion tensor $\mathbf{D}$ (F), the $yy$ component of the diffusion tensor $\mathbf{D}$ (G) and the $zz$ component of the diffusion tensor $\mathbf{D}$ (H) are shown. In the third row,  the tumor tissues segmentation into a viable (brown color), a necrotic component (black color) and edema (pink color) (I), the coefficient of the first spherical harmonics decomposition of the diffusion tensor, in RGB color scale (J), the coefficient of the second spherical harmonics decomposition of the diffusion tensor, in RGB color scale (K) and the coefficient of the third spherical harmonics decomposition of the diffusion tensor, in RGB color scale (L), are shown. We note that higher order spherical harmonics coefficients indicate regions of higher anisotropy, which correspond to the white matter region in subfigures (K) and (L). 

\subsection{Test case 1 - Before surgery}
\begin{figure}[ht!]
\includegraphics[width=1.0\linewidth]
{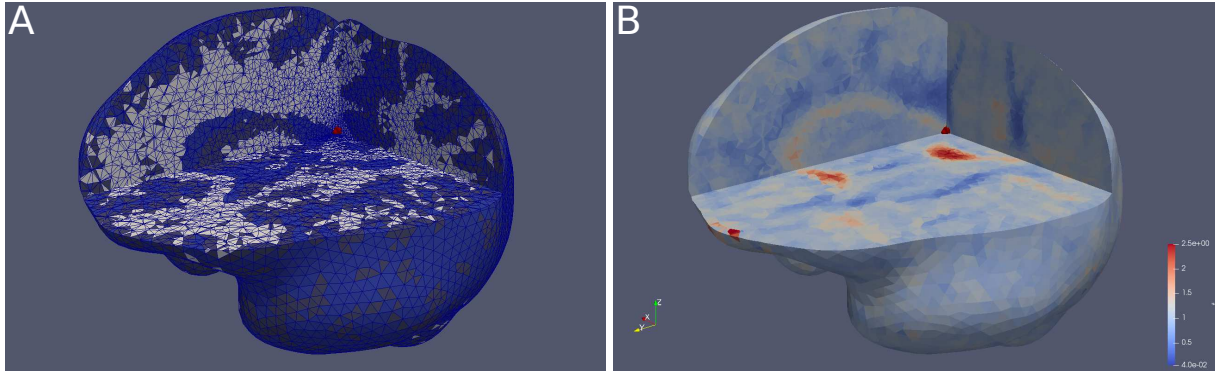}
\centering
\caption{Initial viable cells distribution (brown color) together with the brain computational mesh and the brain tissue labels (white color for WM, gray color for GM, blue color for CSF) (A); initial viable cells distribution (brown color) together with the $xx$ component of the tensor $\mathbf{T}$ (B).}
\label{fig:2}
\end{figure}

In Figure \ref{fig:2}--A we show the initial viable cells distribution (in brown color) for \textbf{Test case 1}, together with the brain computational mesh and the brain tissue labels (white color for WM, gray color for GM, blue color for CSF); in Figure \ref{fig:2}--B we show the same initial viable cells distribution (in brown color) together with the $xx$ component of the tensor $\mathbf{T}$ (which ranges from $0$ to $3$). We take as initial condition a viable cells spherical distribution of radius $2.5$ mm placed in the corpus callosum behind the ventricles (in anterior--posterior orientation), with concentration value $\phi_{v}(0)=0.6$ inside the sphere and zero outside, while we take $\phi_{d}(0)=0$, $\phi_{a}(0)=0$, $n(0)=1$, $c(0)=0$ everywhere. The value $\phi_{v}(0)=0.6$ is higher than the equilibrium value of a uniform state for the viable cells phase with $\phi_{d}(0)=0$, $\phi_{a}(0)=0$, $c(0)=0$ and $n(0)=\delta_n=0.33$, which can be calculated to be $\sim 0.54$ for \textbf{Case 1} and $\sim 0.19$ for \textbf{Case 2}. This means that, both for \textbf{Case 1} and \textbf{Case 2}, the nutrient will be consumed strongly by the viable cells pushing the nutrient below the hypoxia threshold throughout the dynamics. 

\begin{figure}[ht!]
\includegraphics[width=1.0\linewidth]
{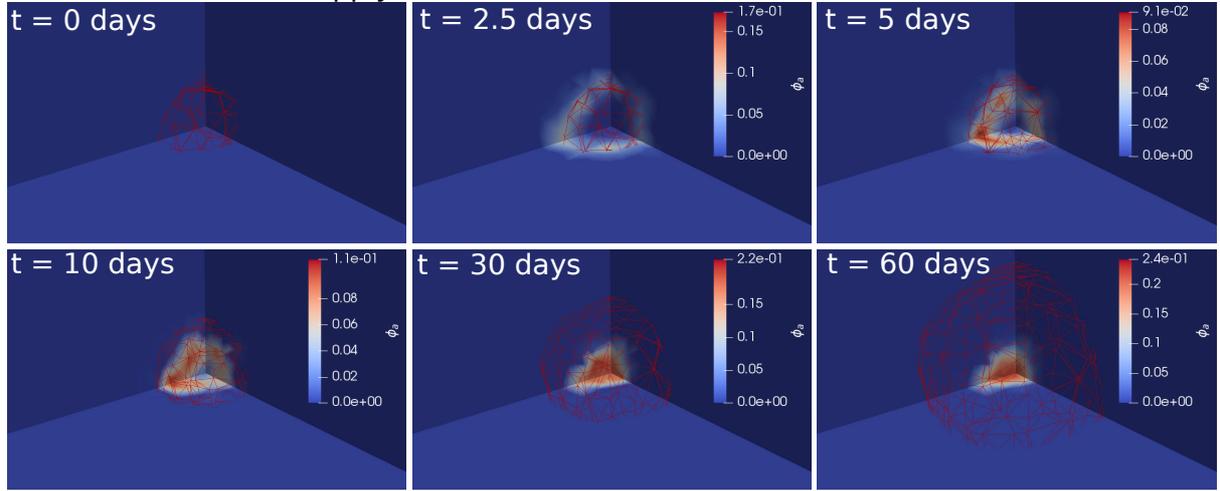}
\centering
\caption{Values of $\phi_a$, superposed to the wireframe plot of the isosurface $\phi_v=0.1$ (in brown color), at time instants $t=0, 2.5, 5, 10, 30, 60$ days, for both the set of parameters of Case 1 and Case 2.}
\label{fig:3}
\end{figure}

In Figure \ref{fig:3} we plot the endothelial cells distribution $\phi_a$, together with the boundary surface of the viable cells distribution $\phi_v$ (in brown color), at different time instants along the system evolution, from $t=0$ to $t=60$ days. The boundary of the viable cells phase is defined here as the isosurface $\phi_v=0.1$. The results for both \textbf{Case 1}, corresponding to a high nutrient supply inside the tumor tissues, and  \textbf{Case 2}, corresponding to a low nutrient supply inside the tumor tissues, are reported.
We observe that, both in \textbf{Case 1} and \textbf{Case 2}, new tumor induced vasculature is produced near the tumor boundary at the earlier stages of evolution. At later time, while the tumor is expanding, the new vasculature moves following the chemotactic gradient of the angiogenetic factor towards the tumor core, infiltrating the tumor tissues. In \textbf{Case 1}, the concentration of the endothelial cells is of the order of magnitude of $0.01$, being the angiogenetic process negligible in this case, while in \textbf{Case 2} the concentration of new vasculature infiltrated in the tumor tissue is of the order of magnitude of $0.1$, i.e. one order of magnitude more concentrated then normal vasculature in the healthy tissue. In this latter case, the angiogenetic process is effective in contrasting the hypoxia regime of viable tumor cells in the tumor core, and the newly induced nutrient supply causes a wider expansion and tumor infiltration in the surrounding tissues.

\begin{figure}[ht!]
\includegraphics[width=1.0\linewidth]
{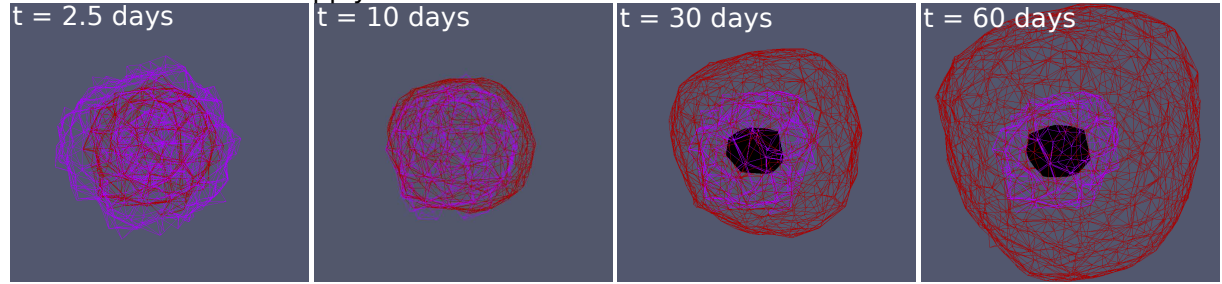}
\centering
\caption{Plot of the isosurfaces $\phi_v=0.1$ (brown color), $\phi_d=0.1$ (black color), $\phi_a=0.02$ (lilac color), at time instants $t=2.5, 10, 30, 60$ days, for both the set of parameters of Case 1 and Case 2. }
\label{fig:4}
\end{figure}

In Figure \ref{fig:4} we plot the isosurfaces $\phi_v=0.1$, $\phi_d=0.1$, $\phi_a=0.02$, identifying the boundary of the viable and necrotic tumor region and the extension of the region interested by the angiogenetic process (considering that $\bar{\phi}_a=0.04$), at different time instants and both for \textbf{Case 1} and \textbf{Case 2}.
We observe that in \textbf{Case 2}, corresponding to the situation with low nutrient supply inside the tumor tissues, a necrotic core develops at later instant of times. Moreover, the additional nutrient supply coming from the tumor induced vasculature causes a wider expansion of the viable phase.

\begin{figure}[ht!]
\includegraphics[width=1.0\linewidth]
{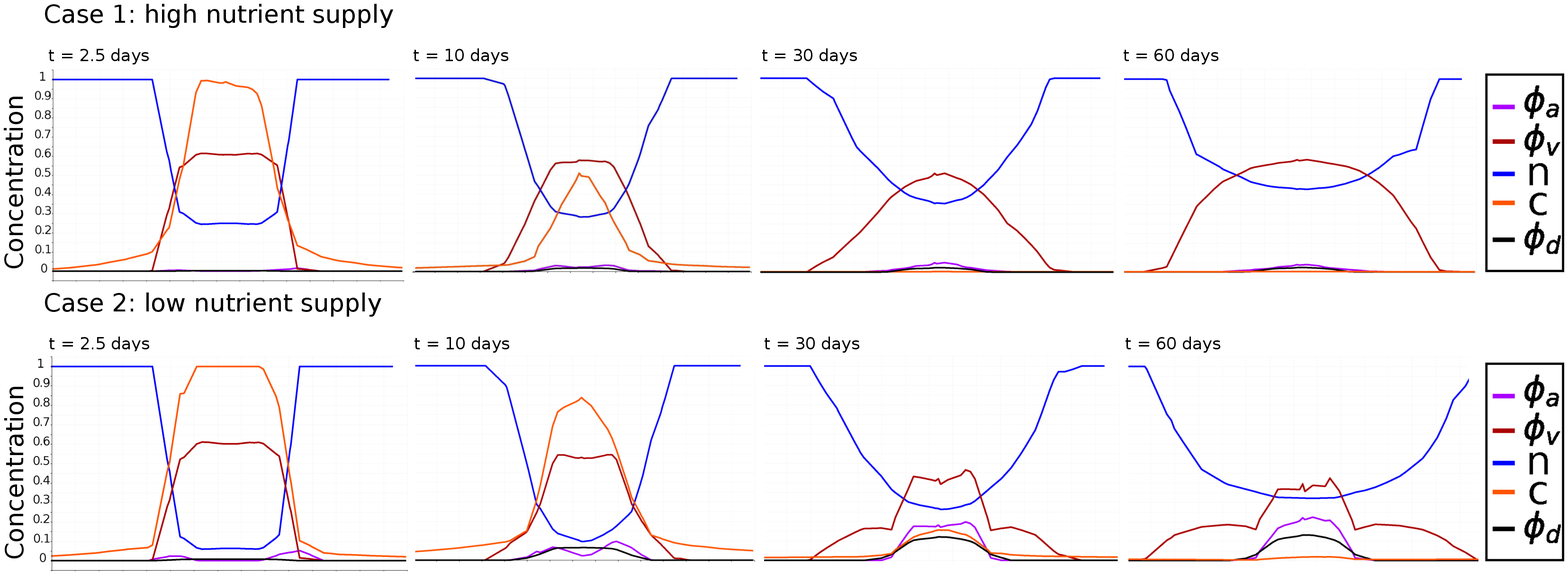}
\centering
\caption{Line plots, along a diameter of the initial spherical tumor distribution, of the concentrations $\phi_v, \phi_d, \phi_a, n, c$, at time instants $t=2.5, 10, 30, 60$ days, for both the set of parameters of Case 1 and Case 2.}
\label{fig:5}
\end{figure}

Finally, in Figure \ref{fig:5} we show the line plots, along a diameter of the initial spherical tumor distribution, of the concentrations $\phi_v, \phi_d, \phi_a, n, c$, at different time instants and both for \textbf{Case 1} and \textbf{Case 2}. 
We observe the negligible production of necrotic tumor cells and tumor induced vasculature in \textbf{Case 1}, while in \textbf{Case 2} the endothelial cells, starting from the boundary of the viable cells and in correspondence with the production of the angiogenetic factor, infiltrate the tumor core and cause an increase of the nutrient supply, thus raising the concentration of the viable cells in the tumor core and widening the region of the tumor boundary (raising the level of tumor infiltration in the surrounding tissues).

\subsection{Test case 2 - Before surgery}
\begin{figure}[ht!]
\includegraphics[width=1.0\linewidth]
{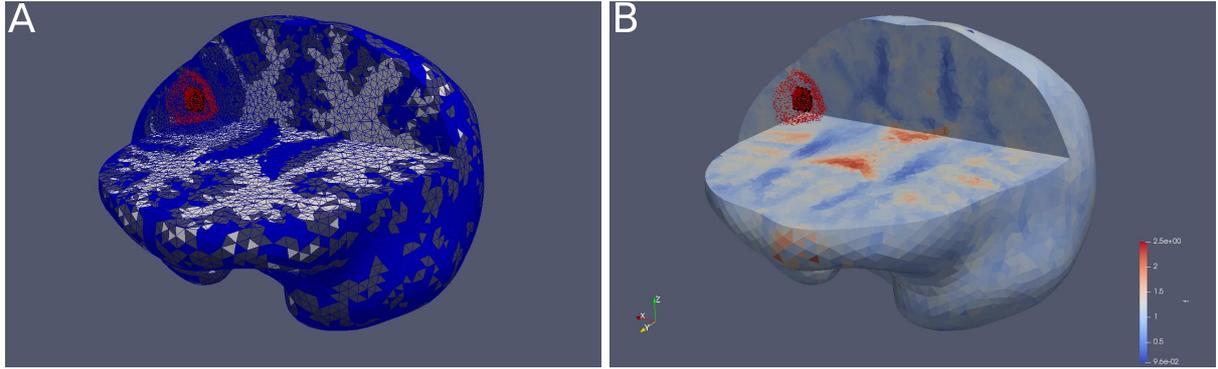}
\centering
\caption{Initial viable cells (brown color) and necrotic cells (black color) distributions, together with the brain computational mesh and the brain tissue labels (white color for WM, gray color for GM, blue color for CSF) (A); initial viable cells (brown color) and necrotic cells (black color) distributions, together with the $xx$ component of the tensor $\mathbf{T}$ (B).}
\label{fig:6}
\end{figure}

In Figure \ref{fig:6}--A we show the initial viable cells distribution (in brown color) and necrotic cells distribution (black color) for \textbf{Test Case 2}, together with the brain computational mesh and the brain tissue labels (white color for WM, gray color for GM, blue color for CSF); in Figure \ref{fig:6}--B we show the same initial viable cells distribution (in brown color) and necrotic cells distribution (black color), together with the $xx$ component of the tensor $\mathbf{T}$ (which ranges from $0$ to $3$). We take as initial condition the viable and necrotic tumor extensions as delineated through the automatic segmentation process described in section \ref{sec:preprocessing}. We set $\phi_{v}(0)=0.6$ inside the viable cells region and zero outside, while we take $\phi_{d}(0)=\bar{\phi}=0.389$ (the equilibrium value of the cells potential) inside the necrotic cells region and zero outside, $\phi_{a}(0)=0$, $n(0)=1$, $c(0)=0$ everywhere. 

\begin{figure}[ht!]
\includegraphics[width=1.0\linewidth]
{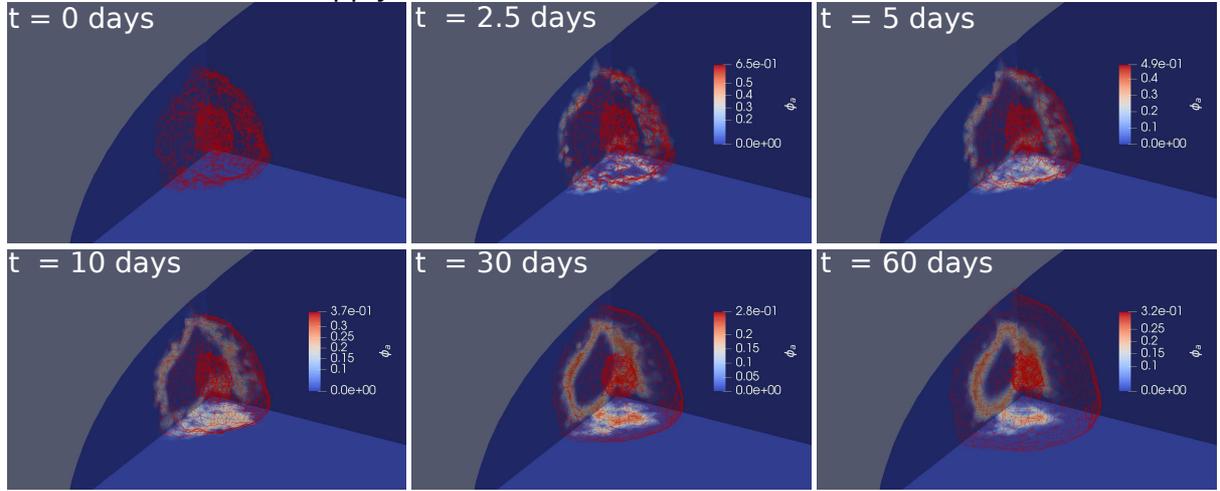}
\centering
\caption{Values of $\phi_a$, superposed to the wireframe plot of the isosurface $\phi_v=0.1$ (in brown color), at time instants $t=0, 2.5, 5, 10, 30, 60$ days, for both the set of parameters of Case 1 and Case 2.}
\label{fig:7}
\end{figure}

In Figure \ref{fig:7} we plot the endothelial cells distribution $\phi_a$, together with the boundary surface of the viable cells distribution $\phi_v$ (in brown color), at different time instants along the system evolution, from $t=0$ to $t=60$ days. The boundary of the viable cells phase is defined here as the isosurface $\phi_v=0.1$. The results for both \textbf{Case 1}, corresponding to a high nutrient supply inside the tumor tissues, and  \textbf{Case 2}, corresponding to a low nutrient supply inside the tumor tissues, are reported.
Similarly to the results reported in Figure \ref{fig:3} for Test Case 1, we observe that, both in \textbf{Case 1} and \textbf{Case 2}, new tumor induced vasculature is produced near the tumor boundary at the earlier stages of evolution, and infiltrates in the tumor tissues by chemotaxis at later times. In \textbf{Case 1}, the concentration of the endothelial cells is one order of magnitude smaller than the same concentration obtained in \textbf{Case 2}. In this latter case, the angiogenetic process compete with the nutrient consumption by the viable cells and it is effective to overpass the hypoxia regime of viable tumor cells in the tumor core.

\begin{figure}[ht!]
\includegraphics[width=1.0\linewidth]
{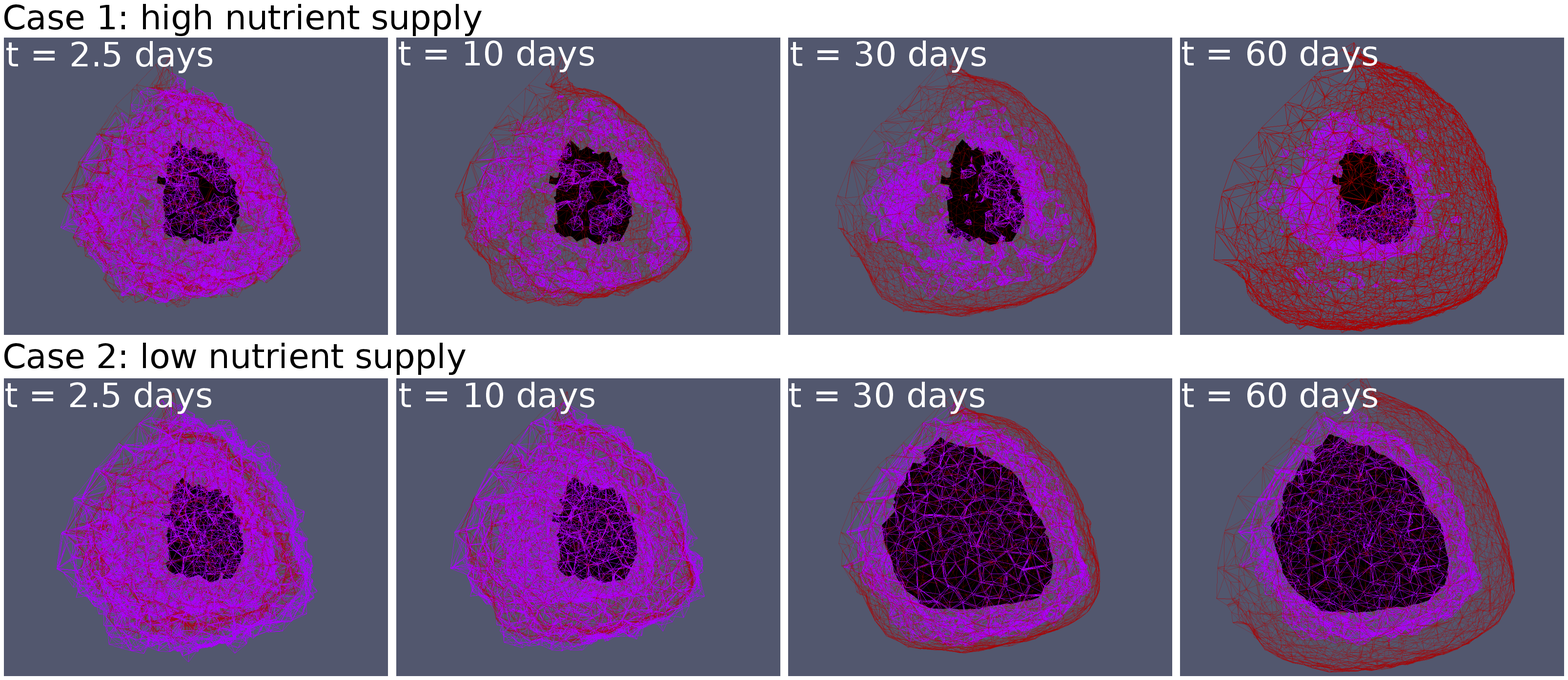}
\centering
\caption{Plot of the isosurfaces $\phi_v=0.1$ (brown color), $\phi_d=0.1$ (black color), $\phi_a=0.02$ (lilac color), at time instants $t=2.5, 10, 30, 60$ days, for both the set of parameters of Case 1 and Case 2. }
\label{fig:8}
\end{figure}

In Figure \ref{fig:8} we plot the isosurfaces $\phi_v=0.1$, $\phi_d=0.1$, $\phi_a=0.02$, identifying the boundary of the viable and necrotic tumor region and the extension of the region interested by the angiogenetic process, at different time instants and both for \textbf{Case 1} and \textbf{Case 2}.
We observe that in \textbf{Case 2}, corresponding to the situation with low nutrient supply inside the tumor tissues, the necrotic core evolves from its initial configuration and expands in a region filled with tumor induced vasculature. 

\begin{figure}[ht!]
\includegraphics[width=1.0\linewidth]
{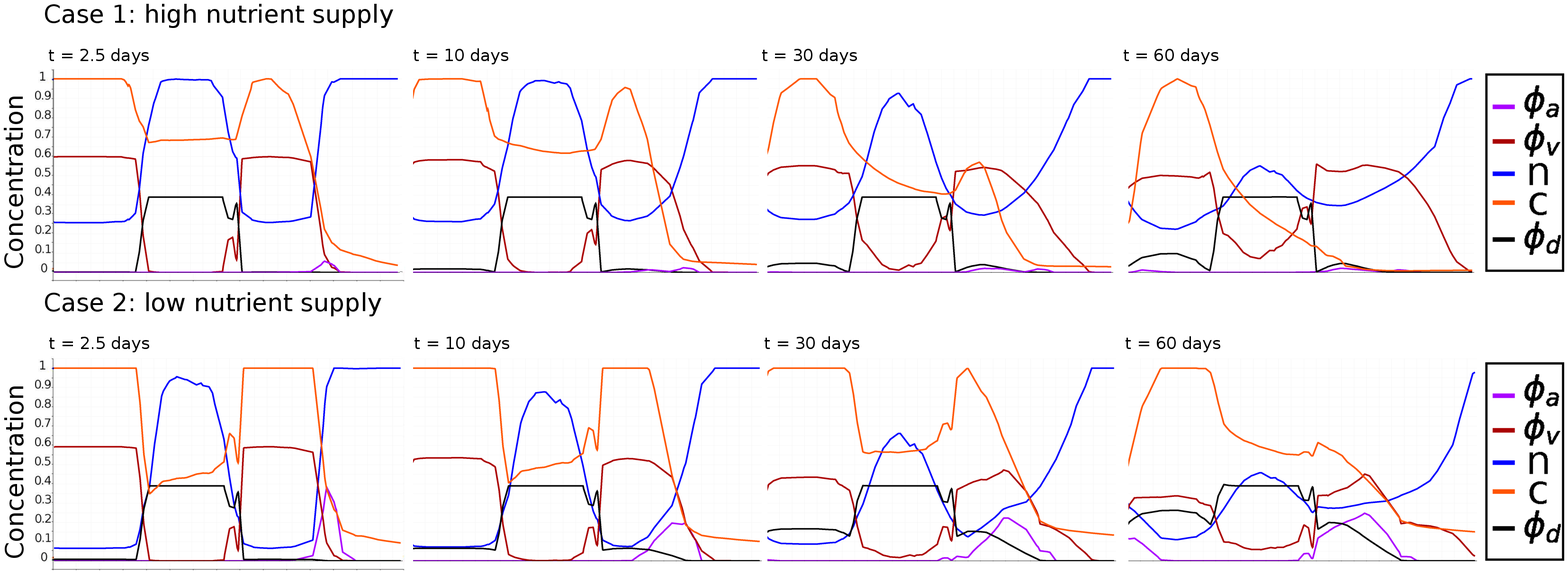}
\centering
\caption{Line plots, along a center--line in the longitudinal extension of the initial tumor distribution (from brain boundary to the brain interior), of the concentrations $\phi_v, \phi_d, \phi_a, n, c$, at time instants $t=2.5, 10, 30, 60$ days, for both the set of parameters of Case 1 and Case 2.}
\label{fig:9}
\end{figure}

Finally, in Figure \ref{fig:9} we show the line plots, along a center--line in the longitudinal extension of the initial tumor distribution (from brain boundary to the brain interior), of the concentrations $\phi_v, \phi_d, \phi_a, n, c$, at different time instants and both for \textbf{Case 1} and \textbf{Case 2}. 
As in Test Case 1, we observe the negligible production of tumor induced vasculature in \textbf{Case 1}, with a fixed necrotic core and an expanding viable cells region. In \textbf{Case 2} the tumor induced vasculature develops at the boundary of the viable cells in the inner region of the brain, (indeed, no healthy tissue with normal vasculature is available at the brain boundary), and in correspondence with the production of the angiogenetic factor. The angiogenetic phase infiltrate the tumor core, characterized by expanding viable cells and a growing necrotic core due to hypoxia, and causes an increase of the nutrient supply, raising the concentration of the viable cells at the tumor boundary and widening the viable cells boundary extension.

\subsection{Test Case 3 - After surgery}
In this test case we study the evolution after surgery of residual tumor particles infiltrated in the boundary of the resection area, and compare the simulation results with the MRI end DSC--MRI data acquired seven months after surgery. During this time span, therapy was administered to the patient according to the standard Stupp protocol. We observe that, from visual inspection of the structural MRI data acquired after surgery, no residual tumor distribution was identified by expert operators (A. G. L., S. L.), while the automatic segmentation tool, described in Section \ref{sec:preprocessing}, identified uncorrectly the resection area visible as an hyperintense region in the postcontrast T1--weighted image as enhancing tumor. We thus identify the tumor distribution after surgery by randomly placing viable tumor cells, with a concentration corresponding to the mechanical equilibrium value $\phi_v(0)=0.389$, on the boundary of the resection area, defined as the boundary of the segmentation map generated by the automatic segmentation tool.

The parameters of the model are chosen as in \textbf{Case 2} in the previous test cases. Further, in order to model the resection area as a wounded region without the presence of normal vasculature, we introduce the (smooth) indicator function $I_R^C$ of the region surrounding the resection area, and redefine the source terms for the nutrient and the new vasculature as 
\begin{align*}
&S_n=\left(V_nI_R^C\left(1-H_r(\phi_v+\phi_d)\right)+V_TI_R^CH_r(\phi_v+\phi_d)(1-\phi_v-\phi_d)\right)(\bar{n}-n)+V_{an}\phi_a(\bar{n}-n)-\delta_v\phi_vn,\\
& \frac{\Gamma_a}{\gamma}=\left(1-H_r(\phi_v+\phi_d)\right)\left(I_R^CV_a\left(1-\phi_v-\phi_d-\phi_a\right)\max(0,c-\delta_c)-k_3\phi_a\right),
\end{align*}
i.e. we put the growth rates due to the normal vasculature equal to zero inside the resection area. Moreover, in order to model therapy, we write 
\[
k_{T,1}(t)=k_{T,2}(t)=k_R(t)+k_C(t),
\]
where $k_R(t)$ and $k_C(t)$ are the temporal profiles of cells death rate associated to the radiotherapy and chemotherapy schedules, as defined in \cite{agosti2}.

\begin{figure}[ht!]
\includegraphics[width=1.0\linewidth]
{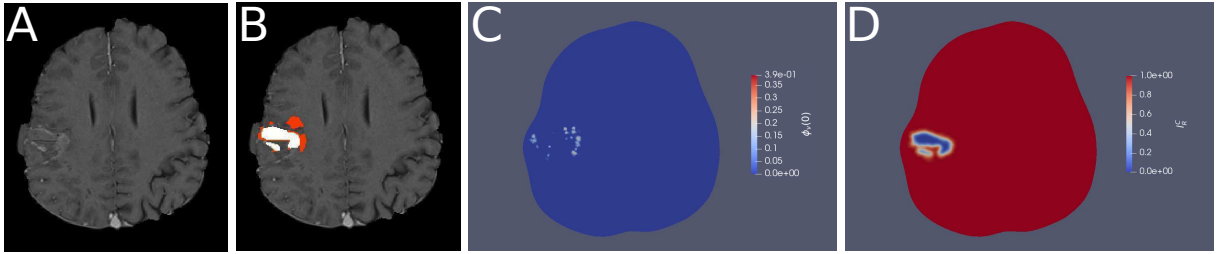}
\centering
\caption{Axial slice from the postcontrast T1-weighted MRI after surgery (A); automatic segmentation of the the tumor tissues into an enhancing component (white color) and a non--enhancing component (red color) (B); initial viable cells distribution $\phi_v(0)$ (C); profile of the indicator function $I_R^C$ (D).}
\label{fig:10}
\end{figure}

In Figure \ref{fig:10}--A we show a representative axial slice of the postcontrast T1-weighted MRI after surgery, while in Figure \ref{fig:10}--B we show the corresponding automatic segmentation of the tumor tissues into an enhancing component (white color) and a non--enhancing component (red color). In Figure \ref{fig:10}--C we show the initial viable cells distribution $\phi_v(0)$ in the same representative axial slice, obtained by randomly placing viable tumor cells at their mechanical equilibrium value in the boundary of the resection area, which gives the initial condition for \textbf{Text Case 3}, together with $\phi_{d}(0)=0$, $\phi_{a}(0)=0$, $c(0)=0$ everywhere, and $n(0)=I_R^C$. In Figure \ref{fig:10}--D we show the profile of $I_R^C$.

\begin{figure}[ht!]
\includegraphics[width=1.0\linewidth]
{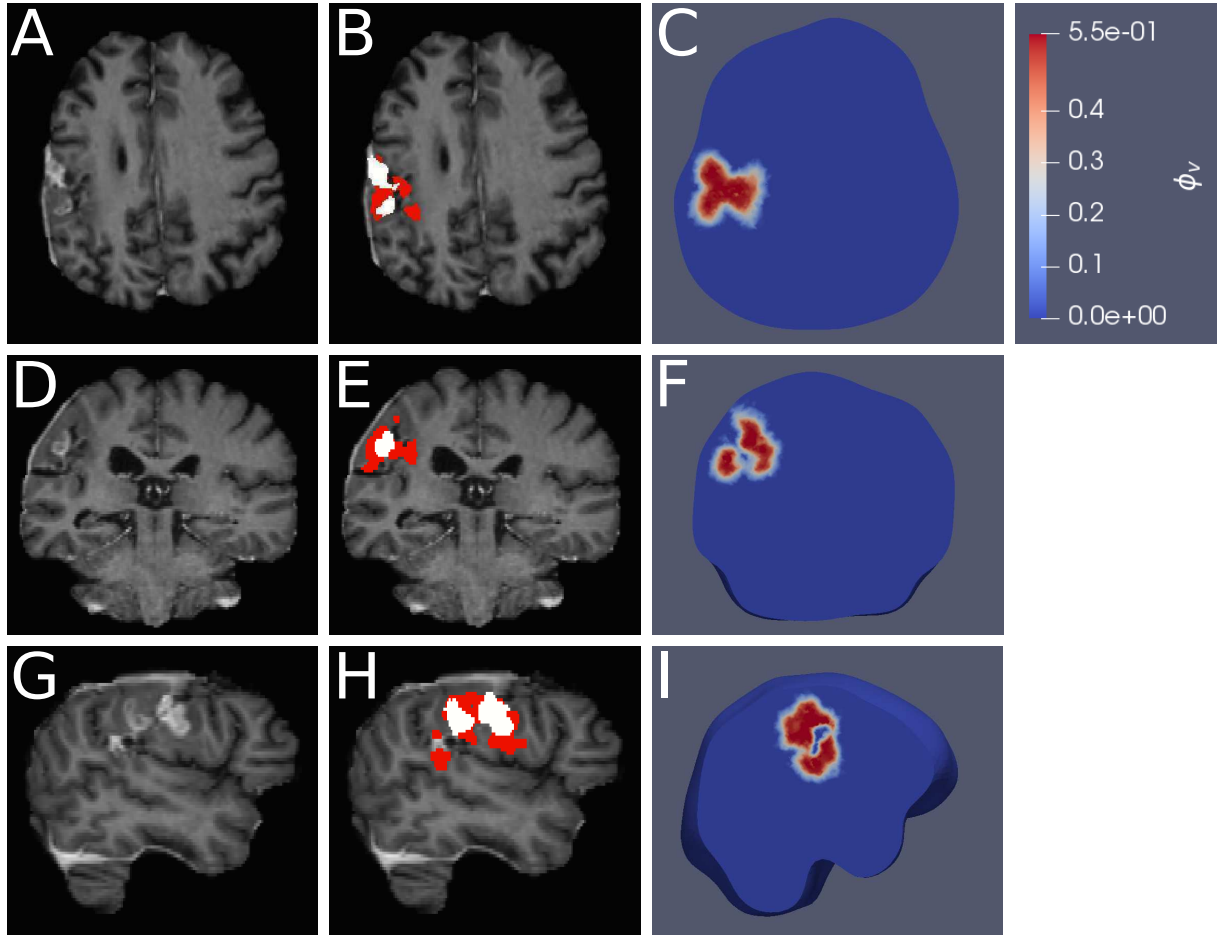}
\centering
\caption{Representative axial (A), coronal (D) and sagittal (G) slices from the postcontrast T1-weighted MRI after therapy; automatic segmentation of the tumor tissues into an enhancing component (white color) and a non--enhancing component (red color) (B, E, H) in the corresponding slices; viable cells distribution $\phi_v$ obtained from numerical simulations (C, F, I) in the corresponding slices.}
\label{fig:11}
\end{figure}

In Figures \ref{fig:11}--A, \ref{fig:11}--D and \ref{fig:11}--G  we show representative axial, coronal and sagittal slices of the postcontrast T1-weighted MRI after therapy, acquired seven months after surgery, together with the corresponding automatic segmentation of the the tumor tissues into an enhancing component (white color) and a non--enhancing component (red color) in Figures \ref{fig:11}--B, \ref{fig:11}--E and \ref{fig:11}--H. We observe that no necrotic component was identified by the automatic segmentation tool in this temporal stage. We qualitatively compare the tumor distributions identified by the automatic segmentation tools with the results obtained by the numerical simulations of our model, shown in the same slices in Figures \ref{fig:11}--C, \ref{fig:11}--F and \ref{fig:11}--I.

We note that the simulation results show the absence of a necrotic component, in accordance with the automatic segmentation results of the MRI data. By comparing the axial slices 11--B and 11--C between data and simulations, we observe that the simulation results reproduce the overall displacement of the tumor distribution identified from data by the automatic segmentation tool, with the simulation results expressing a wider infiltration area of the viable tumor cells in the boundary of the resection area. 
By comparing the coronal slices 11--E and 11--F, we further observe that the numerical simulations predict a greater involvement of tumor extension in the Frontoparietal Operculum, in an hypointense region in the postcontrast T1-weighted image interested by the tumor infiltration process which is not recognized as a tumor region by the automatic segmentation tool. 
Finally, by comparing the sagittal slices 11--H and 11--I, we still observe a greater involvement of tumor extension in the Frontoparietal Operculum predicted by the numerical simulations with respect to the tumor extension identified by the automatic segmentation tool. We also observe that the automatic segmentation tool recognizes an hyperintense region at the inferior portion of the Supramarginal Gyrus as enhancing tumor, which is associated to an independent secondary tumor grown in a separated region with respect to the resection area. Since this secondary tumor is independent from the primary tumor growing at the boundary of the resection area, the tumor distribution from numerical simulations does not cover this region. 

All these considerations suggest that the model predictions are able to correctly identify regions of possible tumor recurrence, with the limitation of not being able to reproduce the onset of satellite tumors in separated regions from the primary tumor or the resection area identified by automatic or manual segmentation at the initial time.
These features paves the way to the use of our computational platform to enhance the diagnostic capabilities of existing automatic segmentation tools based on machine learning methods, as well as to predict tumor recurrence during therapy and optimize the patient--specific therapy schedules. We note from the comparison between Figures 11--H and 11--I that the orientation of the resection area is different between data and simulations; this is due to the fact that the resection area in simulations maintain the initial orientation described from data after surgery, since we do not include in our model the effect of tissue elastic displacements during the tumor growth evolution. 

\begin{figure}[ht!]
\includegraphics[width=1.0\linewidth]
{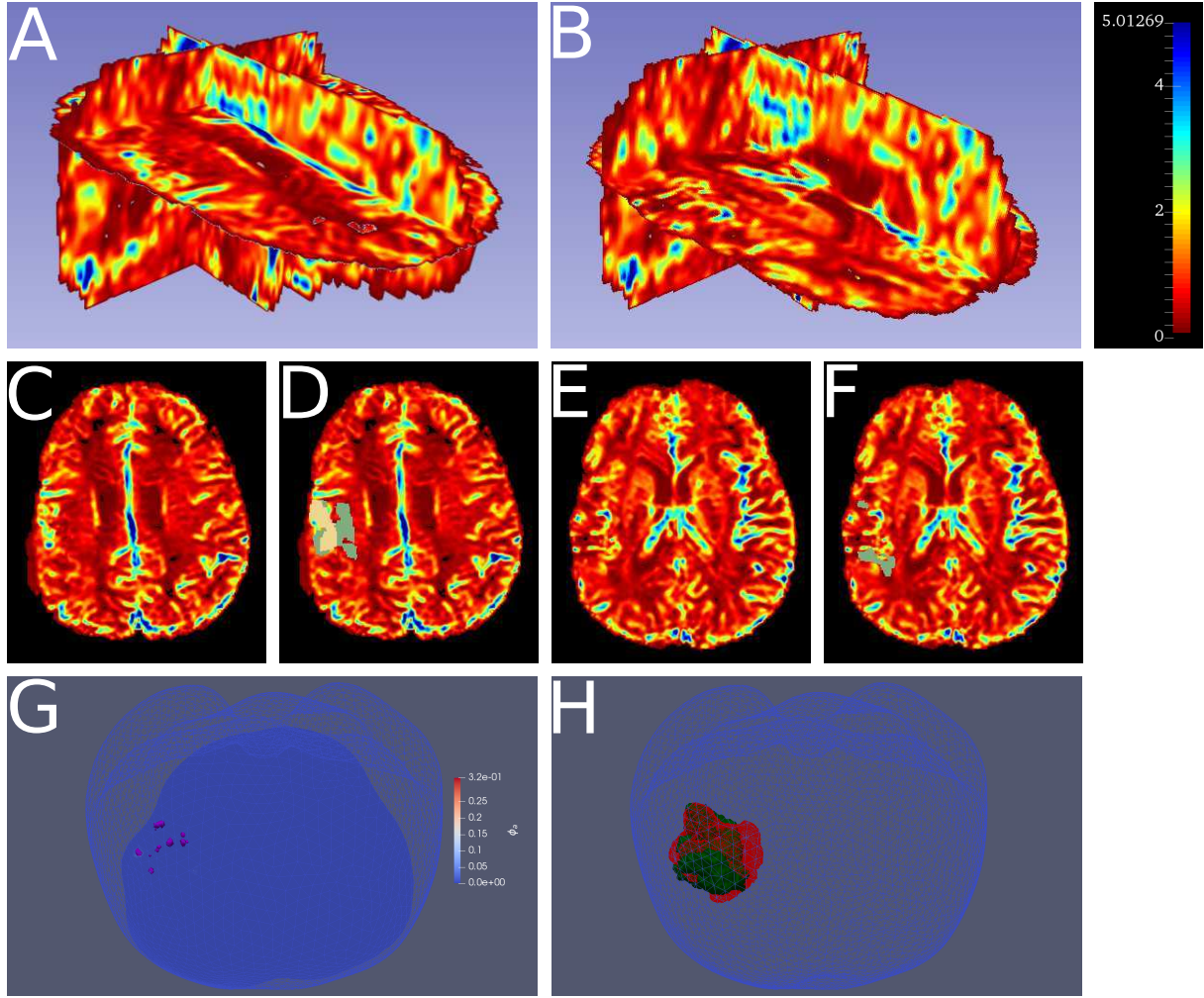}
\centering
\caption{3--d view of two representative axial slices, along the vertical direction, of the rCBV map after therapy (A, B); plane view of the same axial slices reporting the rCBV map (C, E), and reporting the rCBV map with superposed automatic tumor segmentation, with an enhancing component (white color) and a non--enhancing component (green color) (D, F); isosurfaces $\phi_a=0.02$ from the numerical simulations after seven months of evolution, in lilac color, with the upper axial slice highlighted in dark blue color (G); tumor boundary surface from the automatic segmentation map (green color) and isosurface $\phi_v=0.1$ (red color) from the numerical simulations (H).}
\label{fig:12}
\end{figure}
Finally, in Figures \ref{fig:12}--A, \ref{fig:12}--B,  \ref{fig:12}--C and \ref{fig:12}--E  we show a 3--d view and a plane view of two representative axial slices, along the vertical direction, of the rCBV map, calculated from the DSC--MRI data acquired after therapy as described in Section \ref{sec:preprocessing}. We also show the same axial slices reporting the rCBV map with superposed automatic tumor segmentation in Figures \ref{fig:12}--D and \ref{fig:12}--F, with an enhancing component (white color) and a non--enhancing component (green color). In Figure \ref{fig:12}--G we plot the isosurfaces $\phi_a=0.02$, in lilac color, defining the boundaries of the regions occupied by the tumor induced vasculature, obtained from the numerical simulations after seven months of evolution, highlighting also the axial slice corresponding to the one reported in Figure \ref{fig:12}--C in dark blue color for ease of comparison with the rCBV maps. In Figure \ref{fig:12}--H we finally show the tumor boundary surface from the automatic segmentation map (green color) and the isosurface $\phi_v=0.1$ (red color), corresponding to the tumor boundary obtained from the numerical simulations.

By comparing the rCBV maps on the two reported representative axial slices and the numerical solutions for the distribution of the tumor induced vasculature, we observe that the numerical simulations predict the presence of an effective angiogenetic activity near the boundary of the tumor distribution in the same regions which highlight CBV enhancement in the peritumoral area from the perfusion data. Also, we observe that the rCBV peak value inside the angiogenetic regions is one order of magnitude greater than the average rCBV value in the healthy white matter tissue, which is correctly reproduced by the numerical simulations. Indeed, the peak value of $\phi_a$ in the dark blue slice reported in Figure  \ref{fig:12}--G, as shown by the color scale bar, is one order of magnitude greater than $\bar{\phi}_a$. Hence we conclude that the numerical simulations of the proposed model correctly predict both the extension of the angiogenetic process and its intensity in terms of vessels concentration. Finally, in Figure  \ref{fig:12}--H we observe a good superposition between the tumor boundary surface from the automatic segmentation map (green color) and the viable cells extension (red color) from the numerical simulations. A quantitative comparison between simulations and data will be investigated in a forthcoming paper, including model parameters optimization and exploiting more clinical insights.

\section{Conclusions}
In this paper we proposed a new continuous multiphase mixture model describing tumor growth with angiogenesis, including a viable and a necrotic tumor phase, a liquid phase and an angiogenetic phase, coupled with the evolution of a massless nutrient and an angiogenetic factor. The model takes into account the cells--cells and cells--matrix adhesion properties and the infiltrative mechanics of tumor--induced vasculature in the tumor tissues, guided by tumor and endothelial cells proliferation and by the chemotaxis, haptotaxis and porous filtration phenomena with an underlying anisotropy of the tumor and healty tissues, determined by fiber orientations in the tissues. The model is informed by neuroimaging data (MRI, DTI and DSCE--MRI perfusion data) giving informations about the patient--specific geometry and the tissue microstructure, the fiber distributions, the distribution of the different tumor components and the distribution of tumor induced vasculature. We also proposed and implemented a finite element approximation of the model equations, based on a system of coupled discrete variational inequalities, which preserves the qualitative properties of the continuous solutions, in order to run numerical simulations on patient--specific test cases of a patient affected by GBM.  The study of such a mathematical and computational model, which aims to predict the patient--specific tumor tissue oxygenation and tumor--induced microvasculature, together with the tumor infiltration in the peritumoral area, is of crucial importance to assess the patient response to both systemic and radiation therapy and to design optimized therapeutic schedules, especially for GBM, which is highly vascular and infiltrative and for which the application of standard therapy schedules is invariably followed by tumor recurrence. 

Our numerical simulations, in accordance with observations from both in--vitro and in--vivo experiments, showed that the viable tumor cells infiltrates the surrounding tissues following the preferential directions of white matter fiber tracts, while new tumor induced vasculature is produced near the tumor boundary and, following the chemotactic gradient of the angiogenetic factor towards the tumor core, infiltrates the tumor tissues. In situations with low nutrient availability inside the tumor core, a necrotic core develops and the angiogenetic process is effective in contrasting the hypoxia regime of viable tumor cells, with the newly induced nutrient supply causing a wider tumor infiltration in the surrounding tissues. Investigating a test case after surgery and comparing the simulation results with neuroimaging data acquired after therapy, we showed that our model correctly predicts the overall extension of the tumor and tumor--induced vasculature distributions and the intensity
of the angiogenetic process, with the simulation results expressing a wider infiltration area of the viable tumor cells in the boundary of the resection area with respect to the region recognized as tumor by a state of the art automatic segmentation tool.
These features paves the way to the use of our computational platform to enhance the diagnostic capabilities of existing automatic segmentation tools based on machine learning methods, as well as to assist the clinicians in predicting tumor recurrence during therapy, in assessing the therapy outcomes and in designing optimal patient--specific therapeutic schedules.

The proposed model has some limitations, which should be addressed in future developments. Firstly, the parameters in the model were chosen in a range of biological parameters reported in literature, with no empirical estimates available for the parameters in the TAF equation and for the parameters giving the nutrient supply from the tumor induced vasculature. Secondly, the model is not able to predict the onset of secondary independent tumors initiating in regions separated from the ones identified by the initial segmentation process. Finally, the implemented therapy effects in the model are not diversified between the different tumor tissues and depending on the tissue oxygenation and on the tumor and peritumoral microcirculation. Future developments will concern the optimization of the model parameters in a patient--specific manner based on longitudinal neuroimaging data, as well as the inclusion of metastasis dynamics and an implementation of more realistic therapy effects.  

\section{Acknowledgements}
This research has been performed in the framework of the MIUR-PRIN Grant 2020F3NCPX ''Mathematics for industry 4.0 (Math4I4)''. The present paper also benefits from the support of the GNAMPA (Gruppo Nazionale per l'Analisi Matematica, la Probabilit\`a e le loro Applicazioni) of INdAM (Istituto Nazionale di Alta Matematica).
%
\bibliographystyle{plain}
\bibliography{biblio_GBMMP} 

\begin{thebibliography}{10}

\bibitem{agosti1}
A.~Agosti, P.~F. Antonietti, P.~Ciarletta, M.~Grasselli, and M.~Verani.
\newblock A {C}ahn-{H}illiard--type equation with application to tumor growth
  dynamics.
\newblock {\em Math. Methods Appl. Sci.}, 40(18):7598--7626, 2017.
\newblock doi: \url{https://doi.org/10.1002/mma.4548}.

\bibitem{agosti2}
A.~Agosti, C.~Cattaneo, C.~Giverso, D.~Ambrosi, and P.~Ciarletta.
\newblock A computational framework for the personalized clinical treatment of
  glioblastoma multiforme.
\newblock {\em ZAMM}, 98(12):2307--2327, 2018.
\newblock doi: \url{https://doi.org/10.1002/zamm.201700294}.

\bibitem{agosti3}
A.~Agosti, P.~Ciarletta, H.~Garcke, and M.~Hinze.
\newblock Learning patient--specific parameters for a diffuse interface
  glioblastoma model from neuroimaging data.
\newblock {\em Math. Methods Appl. Sci.}, 135(15):8945--8979, 2020.
\newblock doi: \url{https://doi.org/10.1002/mma.6588}.

\bibitem{fenics}
M.~S. Alnaes, J.~Blechta, J.~Hake, A.~Johansson, B.~Kehlet, A.~Logg,
  C.~Richardson, J.~Ring, M.~E. Rognes, and G.~N. Wells.
\newblock The {FE}ni{CS} project version 1.5.
\newblock {\em Archive of Numerical Software}, 3, 2015.
\newblock doi: \url{https://doi.org/10.11588/ans.2015.100.20553}.

\bibitem{intro16}
A.R.A. Anderson and M.A.J. Chaplain.
\newblock Continuous and discrete mathematical models of tumor-induced
  angiogenesis.
\newblock {\em Bull. Math. Biol.}, 60(5):857--899, 1998.
\newblock doi: \url{https://doi.org/10.1006/bulm.1998.0042}.

\bibitem{eddy}
J.~L. Andersson and S.~N. Sotiropoulos.
\newblock An integrated approach to correction for off-resonance effects and
  subject movement in diffusion mr imaging.
\newblock {\em NeuroImage}, 125:1063--1078, 2015.

\bibitem{vmtk}
L.~Antiga, M.~Piccinelli, L.~Botti, B.~Ene-Iordache, A.~Remuzzi, and D.A.
  Steinman.
\newblock An image-based modeling framework for patient-specific computational
  hemodynamics.
\newblock {\em Med. Biol. Eng. Comput.}, 46:1097--1112, 2008.

\bibitem{bv1}
T.~W. Barber, J.~A. Brockway, and L.~S. Higgins.
\newblock The density of tissues in and about the head.
\newblock {\em Acta Neurologica Scandinavica}, 46(1):85--92, 1970.
\newblock doi: \url{https://doi.org/10.1111/j.1600-0404.1970.tb05606.x}.

\bibitem{perfusion2}
J.~L. Boxerman, K.~M. Schmainda, and R.~M. Weisskoffc.
\newblock Relative cerebral blood volume maps corrected for contrast agent
  extravasation significantly correlate with glioma tumor grade, whereas
  uncorrected maps do not.
\newblock {\em AJNR Am. J. Neuroradiol.}, 27(4):859--867, 2006.

\bibitem{intro15}
M.~A.~J. Chaplain and A.~M. Stuart.
\newblock A model mechanism for the chemotactic response of endothelial cells
  to tumour angiogenesis factor.
\newblock {\em IMA J. Math. Appl. Med. Biol.}, 10(3):149--168, 1993.
\newblock doi: \url{https://doi.org/10.1093/imammb/10.3.149}.

\bibitem{chaplain}
M.A.J. Chaplain, S.~M. Giles, B.D. Sleeman, and R.J. Jarvis.
\newblock A mathematical analysis of a model for tumour angiogenesis.
\newblock {\em J. Math. Biol.}, 33(7):744--770, 1995.
\newblock doi: \url{https://doi.org/10.1007/BF00184647}.

\bibitem{onsager3}
C.~Chatelain, T.~Balois, P.~Ciarletta, and M.~Ben~Amar.
\newblock Emergence of microstructural patterns in skin cancer: a phase
  separation analysis in a binary mixture.
\newblock {\em New J. Phys.}, 13:115013, 2011.
\newblock doi: \url{https://doi.org/10.1088/1367-2630/13/11/115013}.

\bibitem{intro11}
S.~Das and P.~A. Marsden.
\newblock Angiogenesis in glioblastoma.
\newblock {\em N. Engl. J. Med.}, 396(16):1561--1563, 2013.
\newblock doi: \url{https://doi.org/10.1056/NEJMcibr1309402}.

\bibitem{intro14}
A.~S. Deakin.
\newblock Model for initial vascular patterns in melanoma transplants.
\newblock {\em Growth}, 40(2):191--201, 1976.

\bibitem{onsager2}
M.~Doi.
\newblock {\em Soft Matter Physics}.
\newblock Oxford University Press, 2013.
\newblock doi: \url{https://doi.org/10.1093/acprof:oso/9780199652952.001.0001}.

\bibitem{temam}
I.~Ekeland and R.~Temam.
\newblock {\em Convex analysis and variational problems}.
\newblock SIAM, 1999.

\bibitem{intro1}
J.~Folkman.
\newblock Tumor angiogenesis.
\newblock {\em Adv. Cancer Res.}, 43:175--203, 1985.
\newblock doi: \url{https://doi.org/10.1016/s0065-230x(08)60946-x}.

\bibitem{intro2}
J.~Folkman and M.~Klagsbrun.
\newblock Angiogenic factors.
\newblock {\em Science}, 235(4787):442--447, 1987.
\newblock doi: \url{https://doi.org/10.1126/science.2432664}.

\bibitem{intro17}
M.~A. Fontelos, A.~Friedman, and B.~Hu.
\newblock Mathematical analysis of a model for the initiation of angiogenesis.
\newblock {\em SIAM J. MATH. ANAL.}, 33(6):1330--1355, 2002.
\newblock doi: \url{https://doi.org/10.1137/S0036141001385046}.

\bibitem{frieboes}
H.~B. Frieboes, F.~Jin, Y.~L. Chuang, S.~M. Wise, J.~S. Lowengrub, and
  V.~Cristini.
\newblock Three-dimensional multispecies nonlinear tumor growth--ii: Tumor
  invasion and angiogenesis.
\newblock {\em J. Theor. Biol.}, 264(4):1254--1278, 2010.
\newblock doi: \url{https://doi.org/10.1016/j.jtbi.2010.02.036}.

\bibitem{intro20}
H.~B. Frieboes, F.~Jin, Y.~L. Chuang, S.~M. Wise, J.~S. Lowengrub, and
  V.~Cristini.
\newblock Three-dimensional multispecies nonlinear tumor growth—ii: Tumor
  invasion and angiogenesis.
\newblock {\em J. Theor. Biol.}, 264(4):1254--1278, 2010.
\newblock doi: \url{https://doi.org/10.1016/j.jtbi.2010.02.036}.

\bibitem{garcke1}
H.~Garcke, K.F. Lam, E.~Sitka, and V.~Styles.
\newblock A {C}ahn--{H}illiard--{D}arcy model for tumour growth with chemotaxis
  and active transport.
\newblock {\em Math Models Methods Appl Sci}, 26(6):1095--1148, 2016.
\newblock doi: \url{https://doi.org/10.1142/S0218202516500263}.

\bibitem{intro7}
R.J. Gillies, P.A. Schornack, T.W. Secomb, and N.~Raghunand.
\newblock Causes and effects of heterogeneous perfusion in tumors.
\newblock {\em Neoplasia}, 1(3):197--207, 1999.
\newblock doi: \url{https://doi.org/10.1038/sj.neo.7900037}.

\bibitem{intro8}
L.H. Gray, A.D. Conger, M.~Ebert, S.~Hornsey, and O.C.A. Scott.
\newblock The concentration of oxygen dissolved in tissues at the time of
  irradiation as a factor in radiotherapy.
\newblock {\em Br. J. Radiol}, 26(312):638--648, 1953.
\newblock doi: \url{https://doi.org/10.1259/0007-1285-26-312-638}.

\bibitem{gurtin1}
M.~E. Gurtin.
\newblock Generalized {G}inzburg--{L}andau and {C}ahn--{H}illiard equations
  based on a microforce balance.
\newblock {\em Phys. D}, 92(3--4).

\bibitem{bv2}
P.~J. Harrison, N.~Freemantle, and J.~R. Geddes.
\newblock Meta-analysis of brain weight in schizophrenia.
\newblock {\em Schizophrenia Research}, 64(1):25--34, 2003.
\newblock doi: \url{https://doi.org/10.1016/s0920-9964(02)00502-9}.

\bibitem{intro18}
T.~A.~M. Heck, M.~M. Vaeyens, and H.~Van~Oosterwyck.
\newblock Computational models of sprouting angiogenesis and cell migration:
  Towards multiscale mechanochemical models of angiogenesis.
\newblock {\em Math. Model. Nat. Phenom.}, 10(1):108--141, 2015.
\newblock doi: \url{https://doi.org/10.1051/mmnp/201510106}.

\bibitem{intro9}
D.~A.~II Hormuth, C.M. Phillips, C.~Wu, E.~A. B.~F. Lima, G.~Lorenzo, P.~K.
  Jha, A.~M. Jarrett, J.~T. Oden, and T.~E. Yankeelov.
\newblock Biologically-based mathematical modeling of tumor vasculature and
  angiogenesis via time-resolved imaging data.
\newblock {\em Cancers}, 13(12):3008, 2021.
\newblock doi: \url{https://doi.org/10.3390/cancers13123008}.

\bibitem{nnunet}
F.~Isensee, P.F. Jaeger, S.A.A. Kohl, J.~Petersen, and K.~H. Maier-Hein.
\newblock nnu-net: a self-configuring method for deep learning-based biomedical
  image segmentation.
\newblock {\em Nat Methods}, 18:203--211, 2021.
\newblock doi: \url{https://doi.org/10.1038/s41592-020-01008-z}.

\bibitem{intro6}
R.K. Jain, E.~di~Tomaso, D.G. Duda, J.S. Loeffler, A.G. Sorensen, and T.T.
  Batchelor.
\newblock Angiogenesis in brain tumours.
\newblock {\em Nat. Rev. Neurosci.}, 8(8):610--622, 2007.
\newblock doi: \url{https://doi.org/10.1038/nrn2175}.

\bibitem{intro13}
B.~S. Jang, S.~H. Jeon, I.~H. Kim, and I.~A. Kim.
\newblock Prediction of pseudoprogression versus progression using machine
  learning algorithm in glioblastoma.
\newblock {\em Sci. Rep.}, 8:12516, 2018.
\newblock doi: \url{https://doi.org/10.1038/s41598-018-31007-2}.

\bibitem{bv3}
K.~K. Kaisti, J.~W. Langsj{\"o}, S.~Aalto, V.~Oikonen, H.~Sipil{\"a},
  M.~Ter{\"a}s, S.~Hinkka, L.~Mets{\"a}honkala, and H.~Scheinin.
\newblock Effects of sevoflurane, propofol, and adjunct nitrous oxide on
  regional cerebral blood flow, oxygen consumption, and blood volume in humans.
\newblock {\em Anesthesiology}, 99(3):603--613, 2003.
\newblock doi: \url{https://doi.org/10.1097/00000542-200309000-00015}.

\bibitem{intro10}
R.F. Kallman and M.J. Dorie.
\newblock Tumor oxygenation and reoxygenation during radiation theraphy: Their
  importance in predicting tumor response.
\newblock {\em Int. J. Radiat. Oncol.}, 12(4):681--685, 1986.
\newblock doi: \url{https://doi.org/10.1016/0360-3016(86)90080-5}.

\bibitem{gibbs}
E.~Kellner, B~Dhital, V.G. Kiselev, and M.~Reisert.
\newblock Gibbs-ringing artifact removal based on local subvoxel-shifts.
\newblock {\em Magnetic Resonance in Medicine}, 76:1574--1581, 2016.
\newblock doi: \url{https://doi.org/10.1002/mrm.26054}.

\bibitem{marchingcubes}
W.E. Lorensen and H.E. Cline.
\newblock Marching cubes: A high resolution 3d surface construction algorithm.
\newblock {\em Computer Graphics}, 21(4):163--169, 1987.

\bibitem{onsager1}
L.~Onsager.
\newblock Reciprocal relations in irreversible processes. i.
\newblock {\em Phys. Rev.}, 37(4).

\bibitem{intro3}
N.~Paweletz and M.~Knierim.
\newblock Tumor--related angiogenesis.
\newblock {\em Crit. Rev. Oncol. Hematol.}, 9(3):197--242, 1989.
\newblock doi: \url{https://doi.org/10.1016/s1040-8428(89)80002-2}.

\bibitem{intro22}
J.~R. Petrella and J.~M. Provenzale.
\newblock Mr perfusion imaging of the brain: techniques and applications.
\newblock {\em AJR Am J Roentgenol.}, 175(1):207--219, 2000.
\newblock doi: \url{https://doi.org/10.2214/ajr.175.1.1750207}.

\bibitem{intro19}
M.~Scianna, C.G. Bell, and L.~Preziosi.
\newblock A review of mathematical models for the formation of vascular
  networks.
\newblock {\em J. Theor. Biol.}, 333:174--209, 2013.
\newblock doi: \url{https://doi.org/10.1016/j.jtbi.2013.04.037}.

\bibitem{intro5}
M.~M. Sholley, G.~P. Ferguson, H.~R. Seibel, J.~L. Montour, and J.~D. Wilson.
\newblock Mechanisms of neovascularization. vascular sprouting can occur
  without proliferation of endothelial cells.
\newblock {\em Lab. Invest.}, 51(6):624--634, 1984.

\bibitem{intro4}
M.~M. Sholley, M.~A. Gimbrone, and R.~S. Cotran.
\newblock Cellular migration and replication in endothelial regeneration.
\newblock {\em Lab. Invest.}, 36(1):18--25, 1977.

\bibitem{tetgen}
H.~Si.
\newblock Tetgen, a delaunay-based quality tetrahedral mesh generator.
\newblock {\em ACM Trans. Math. Software}, 41(2):11:1--11:36, 2015.

\bibitem{nucorrect}
J.G. Sled, A.P. Zijdenbos, and A.C. Evans.
\newblock A nonparametric method for automatic correction of intensity
  nonuniformity in mri data.
\newblock {\em IEEE Trans. Med. Imaging}, 17:87--97, 1998.

\bibitem{bet}
S.~M. Smith.
\newblock Fast robust automated brain extraction.
\newblock {\em Hum. Brain Mapp.}, 17:143--155, 2002.

\bibitem{stupp}
R.~et~al. Stupp.
\newblock Radiotherapy plus concomitant and adjuvant temozolomide for
  glioblastoma.
\newblock {\em N. Engl. J. Med.}, 352:987--996, 2005.
\newblock doi: \url{https://doi.org/10.1056/NEJMoa043330}.

\bibitem{intro21}
S.~Supiot, A.~Lisbona, F.~Paris, D.~Azria, and P.~Fenoglietto.
\newblock "dose-painting": myth or reality?
\newblock {\em Cancer Radiother.}, 14(6--7):554--562, 2010.
\newblock doi: \url{https://doi.org/10.1016/j.canrad.2010.06.005}.

\bibitem{suzuki}
T.~Suzuki.
\newblock {\em Free Energy and Self--Interacting Particles. Progress in
  Nonlinear Differential Equations and Their Applications. Volume 62.}
\newblock Birkh\"{a}user Basel, 2005.
\newblock doi: \url{https://doi.org/10.1007/0-8176-4436-9}.

\bibitem{taubin}
G.~Taubin.
\newblock Curve and surface smoothing without shrinkage.
\newblock {\em Proceedings of IEEE International Conference on Computer Vision,
  Cambridge, MA, USA}, pages 852--857, 1995.

\bibitem{deconv}
J.~D. Tournier, F.~Calamante, D.~G. Gadian, and A.~Connelly.
\newblock Direct estimation of the fiber orientation density function from
  diffusion--weighted mri data using spherical deconvolution.
\newblock {\em NeuroImage}, 23:1176--1185, 2004.
\newblock doi: \url{https://doi.org/10.1016/j.neuroimage.2004.07.037}.

\bibitem{mrtrix}
J.~D. Tournier, R.~Smith, D.~Raffelt, R.~Tabbara, T.~Dhollander, M.~Pietsch,
  D.~Christiaens, B.~Jeurissen, C.~H. Yeh, and A.~Connelly.
\newblock Mrtrix3: A fast, flexible and open software framework for medical
  image processing and visualisation.
\newblock {\em NeuroImage}, 202:116137, 2019.
\newblock doi: \url{https://doi.org/10.1016/j.neuroimage.2019.116137}.

\bibitem{intro12}
B.~R.~J. van Dijken, P.~J. van Laar, M.~Smits, J.~W. Dankbaar, R.~H. Enting,
  and A.~van~der Hoorn.
\newblock Perfusion mri in treatment evaluation of glioblastomas: Clinical
  relevance of current and future techniques.
\newblock {\em J. Magn. Reson. Imaging}, 49(1):11--22, 2019.
\newblock doi: \url{https://doi.org/10.1002/jmri.26306}.

\bibitem{denoise}
J.~Veraart, D.S. Novikov, D.~Christiaens, B.~Ades-aron, J.~Sijbers, and
  E.~Fieremans.
\newblock Denoising of diffusion mri using random matrix theory.
\newblock {\em NeuroImage}, 142:394--406, 2016.
\newblock doi: \url{https://doi.org/10.1016/j.neuroimage.2016.08.016}.

\bibitem{tensor2}
J.~Veraart, J.~Sijbers, S.~Sunaert, A.~Leemans, and B.~Jeurissen.
\newblock Weighted linear least squares estimation of diffusion mri parameters:
  strengths, limitations, and pitfalls.
\newblock {\em NeuroImage}, 81:335--346, 2013.

\bibitem{onsager4}
H.~Wang, T.~Qian, and X.~Xu.
\newblock Onsager's variational principle in active soft matter.
\newblock {\em Soft Matter}, 17(13).

\bibitem{perfusion1}
L.~Willats and F.~Calamante.
\newblock The 39 steps: evading error and deciphering the secrets for accurate
  dynamic susceptibility contrast mri.
\newblock {\em NMR Biomed.}, 26(8):913--931, 2013.
\newblock doi: \url{https://doi.org/10.1002/nbm.2833}.

\bibitem{fsl}
M.W. Woolrich, S.~Jbabdi, B.~Patenaude, M.~Chappell, S.~Makni, T.~Behrens,
  C.~Beckmann, M.~Jenkinson, and S.M. Smith.
\newblock Bayesian analysis of neuroimaging data in fsl.
\newblock {\em NeuroImage}, 45:S173--86, 2009.

\bibitem{fast}
Y.~Zhang, M.~Brady, and S.~Smith.
\newblock Segmentation of brain mr images through a hidden markov random field
  model and the expectation-maximization algorithm.
\newblock {\em IEEE Trans Med Imag}, 20(1):45--57, 2001.

\end{thebibliography}
\end{document}